\documentclass[11pt]{amsart}
\usepackage{amsmath,amsthm,amssymb,amsfonts}
\usepackage{hyperref}
\usepackage[enableskew,vcentermath]{youngtab}
\usepackage[all]{xy}

\newcommand{\cellsize}{11}
\newlength{\cellsz} 
\newsavebox{\cell}
\sbox{\cell}{\begin{picture}(\cellsize,\cellsize)
\put(0,0){\line(1,0){\cellsize}} \put(0,0){\line(0,1){\cellsize}}
\put(\cellsize,0){\line(0,1){\cellsize}}
\put(0,\cellsize){\line(1,0){\cellsize}}
\end{picture}}
\newcommand\cellify[1]{\def\thearg{#1}\def\nothing{}%
\ifx\thearg\nothing \vrule width0pt height\cellsz depth0pt\else
\hbox to 0pt{\usebox{\cell} \hss}\fi%
\vbox to \cellsz{ \vss \hbox to \cellsz{\hss$#1$\hss} \vss}}
\newcommand\tableau[1]{\vcenter{\vbox{\let\\\cr
\baselineskip -16000pt \lineskiplimit 16000pt \lineskip 0pt
\ialign{&\cellify{##}\cr#1\crcr}}}}
\newcommand\tabl[1]{\vtop{\let\\\cr
\baselineskip -16000pt \lineskiplimit 16000pt \lineskip 0pt
\ialign{&\cellify{##}\cr#1\crcr}}}

\newcommand{\A}{\mathbb A}
\newcommand{\al}{\alpha}

\newcommand{\B}{\mathbb B}

\newcommand{\can}{\mathrm{can}}

\newcommand{\core}{\mathcal{C}_{2n}}
\newcommand{\C}{\mathbb{C}}
\newcommand{\Core}{\mathcal{C}_n}
\newcommand{\db}{\Phi}

\newcommand{\Des}{\mathrm{Des}}
\newcommand{\DU}{\mathcal{DU}}
\newcommand{\fin}{{\rm fin}}

\newcommand{\geh}{\mathfrak{g}}
\newcommand{\Gs}{\Gamma_s}
\newcommand{\Gw}{\Gamma_w}
\newcommand{\Gr}{{\rm Gr}}
\newcommand{\h}{\mathfrak{h}}

\newcommand{\idd}{{\rm id}}
\newcommand{\inner}[2]{\langle #1\,,\,#2\rangle}
\newcommand{\ins}{\mathrm{I}}
\newcommand{\insh}{\mathrm{in}}

\newcommand{\Inv}{{\rm Inv}}
\newcommand{\Invj}{\Inv|_{W^i}}
\newcommand{\la}{\lambda}
\newcommand{\La}{\Lambda}
\newcommand{\om}{\omega}
\newcommand{\omt}{\widetilde{\omega}}
\newcommand{\outsh}{\mathrm{out}}
\newcommand{\Q}{\mathbb Q}
\newcommand{\RealRoot}{\Delta_{re}}

\newcommand{\strong}{{\rm strong}}
\newcommand{\tb}{\tilde{\beta}}

\newcommand{\too}{\tilde{\omega}}
\newcommand{\tr}{\mathrm{tr}}
\newcommand{\ts}{{\tilde{S}_{2n}}}
\newcommand{\tS}{\tilde{S}_n}
\newcommand{\tu}{\tilde{u}}
\newcommand{\Tab}{\mathcal{T}}
\newcommand{\UD}{\mathcal{UD}}
\newcommand{\vn}{\varnothing}
\newcommand{\weak}{{\rm weak}}
\newcommand{\wless}{\prec}
\newcommand{\Y}{\mathbb{Y}}
\newcommand{\Z}{\mathbb{Z}}

\newtheorem{thm}{Theorem}
\numberwithin{thm}{section}
\newtheorem{lem}[thm]{Lemma}
\newtheorem{prop}[thm]{Proposition}
\newtheorem{cor}[thm]{Corollary}

\theoremstyle{definition}
\newtheorem{ex}[thm]{Example}
\newtheorem{remark}[thm]{Remark}

\author{Thomas Lam}
\address{Department of Mathematics, Harvard University, Cambridge MA
02138 USA} \email{tfylam@math.harvard.edu}
\thanks{T. L. was partially supported by NSF DMS--0600677.}
\author{Mark Shimozono}
\address{Department of Mathematics, Virginia Polytechnic Institute
and State University, Blacksburg, VA 24061-0123 USA}
\email{mshimo@vt.edu}

\thanks{M. S. was partially supported by NSF DMS--0401012.}

\title{Dual graded graphs for Kac-Moody algebras}
\begin{document}
\begin{abstract}
Motivated by affine Schubert calculus, we construct a family of dual
graded graphs $(\Gs,\Gw)$ for an arbitrary Kac-Moody algebra $\geh$.
The graded graphs have the Weyl group $W$ of $\geh$ as vertex set
and are labeled versions of the strong and weak orders of $W$
respectively.  Using a construction of Lusztig for quivers with an
admissible automorphism, we define folded insertion for a Kac-Moody
algebra and obtain Sagan-Worley shifted insertion from
Robinson-Schensted insertion as a special case. Drawing on work of
Proctor and Stembridge, we analyze the induced subgraphs of
$(\Gs,\Gw)$ which are distributive posets.
\end{abstract}
\maketitle

\section{Introduction}
The Robinson-Schensted correspondence is perhaps the most important
algorithm in algebraic combinatorics. It exhibits a bijection
between permutations and pairs of standard Young tableaux of the
same shape. Stanley investigated the class of differential
posets~\cite{Sta} (also studied in~\cite{Fomin:gen}) and Fomin
studied the more general notion of a dual graded
graph~\cite{Fomin:dual} to formalize local conditions which
guarantee the existence of a Robinson-Schensted style algorithm.

In this article, we construct a family of dual graded graphs
$(\Gs,\Gw)$ associated to each Kac-Moody algebra $\geh$. These
graded graphs have as vertex set the Weyl group $W$ of $\geh$. The
pair $(\Gs,\Gw) = (\Gs(\La),\Gw(K))$ depends on a pair $(\La,K)$
where $\La$ is a dominant integral weight and $K$ is a ``positive
integral'' element of the center $Z(\geh)$. In every case $\Gw$ is
obtained by labeling the left weak order of $W$ and $\Gs$ is
obtained by labeling the strong Bruhat order of $W$.

These labelings are motivated by the Schubert calculus for
homogeneous spaces associated to the Kac-Moody group $G$ with Lie
algebra $\geh$. For $w\in W$, let $\xi^w \in H^*(G/B)$ denote the
cohomology Schubert classes of the flag manifold of $G$. If $\Lambda
= \Lambda_i$ is the $i$-th fundamental weight, then an edge $v
\lessdot w$ in $\Gs(\Lambda)$ is labeled with the coefficient of
$\xi^w$ in the product $\xi^{s_i} \xi^v$, also called a Chevalley
coefficient. When $\geh$ is of affine type and $K = K_\can$ is the
canonical central element, the analogous statement holds
(Proposition~\ref{P:homology}) for $\Gw(K)$ with the homology
Schubert classes $\xi_w \in H_*(\Gr)$ of the affine Grassmannian
corresponding to $\geh$ replacing the cohomology classes. Thus the
combinatorics of these graphs encode computations in Schubert
calculus, and the duality of the graded graphs $(\Gs,\Gw)$ is a
combinatorial skeleton of the duality between cohomology and
homology of homogeneous spaces of $G$.

In the case of the affine Grassmannian, the dual graded graph
structure arises from the pair of dual graded Hopf algebras given by
$H_*(\Gr)$ and $H^*(\Gr)$: one may define the down operator by the
action of the homology class $\xi_{s_0}$ on the Schubert basis of
$H^*(\Gr)$ and the up operator by multiplication by $\xi^{s_i}$ for
any fixed simple reflection $s_i$. It is a general phenomenon that
pairs of dual graded combinatorial Hopf algebras yield dual graded
graphs; we shall pursue this in a separate publication \cite{LS}.

Chains in the graded graphs $(\Gs,\Gw)$, which we call strong and
weak tableaux, are natural generalizations of standard Young
tableaux.  To go one speculative step further, we believe that the
generating functions of an appropriate semistandard notion of strong
and weak tableaux would give polynomials which {\it represent}
certain homology and cohomology Schubert classes, in particular for
homogeneous spaces corresponding to maximal parabolics, generalizing
Schur functions, Schur $Q$-functions and the like.  While this
statement is vague in general, it can be made much more precise when
$\geh$ is of affine type, and has already been achieved in one case.

In the case that $\geh$ is of the affine type $A_{n-1}^{(1)}$ our
construction recovers the dual graded graphs that were implicitly
studied in our joint work with Lapointe and Morse~\cite{LLMS}. The
weak and strong tableaux in~\cite{LLMS} are semistandard
generalizations of the corresponding objects here. In~\cite{LLMS},
an affine insertion algorithm was explicitly constructed for
semistandard weak and strong tableaux, and due to \cite{Lam,LLMS} it
is known that the corresponding generating functions do indeed
represent Schubert classes of the affine Grassmannian of type $A$.
In the limit $n \to \infty$ of the $A_{n-1}^{(1)}$ case, our
construction reproduces Young's lattice, which is the self-dual
graded graph that gives rise to the Robinson-Schensted algorithm.

Having constructed the Kac-Moody dual graded graphs we study two
further aspects of these graphs in detail.

The first aspect is motivated by the relation between the
Robinson-Schensted insertion and Sagan-Worley shifted insertion.
Using Lusztig's construction \cite{L} which associates to each
symmetrizable generalized Cartan matrix $A$, a symmetric generalized
Cartan matrix $B$ equipped with an admissible automorphism $\pi$, we
show that any dual graded graph of the form $(\Gs^A,\Gw^A)$ for
$\geh(A)$ can be realized in terms of one of the form
$(\Gs^B,\Gw^B)$ for $\geh(B)$. In particular, for any affine
algebra, any of the dual graded graphs $(\Gs,\Gw)$ can be realized
using a dual graded graph for a simply-laced affine algebra. In
particular we obtain a Schensted bijection for type $C_n^{(1)}$,
using the insertion algorithm of \cite{LLMS} for type
$A_{2n-1}^{(1)}$. As $n$ goes to infinity, the type $C_n^{(1)}$
insertion converges to Sagan-Worley insertion \cite{Sagan,Worley}.
As a related result, we define a notion of mixed insertion for dual
graded graphs equipped with a pair of automorphisms.  This
generalizes Haiman's variants of Schensted insertion known as
left-right, mixed, and doubly-dual insertion \cite{Ha}.

The second aspect we investigate are the induced subgraphs of the
pair $(\Gs,\Gw)$ which are distributive lattices when considered as
posets. These are precisely the conditions under which one may
describe our strong and weak tableaux by ``filling cells with
numbers" as in a usual standard Young tableau. Here we draw on work
of Proctor \cite{Pro} \cite{Pro2} and Stembridge \cite{Ste}, which
classifies the parabolic quotients of Weyl groups of simple Lie
algebras whose left weak orders (or equivalently Bruhat orders) are
distributive lattices. We sharpen these results slightly to show
that in these cases, the distributivity is compatible with the edge
labels of the graphs $(\Gs,\Gw)$; see Section~\ref{sec:dist}. These
distributive parabolic quotients have also appeared recently in the
geometric work of Thomas and Yong~\cite{TY}. They show that in these
cases one may use the \textit{jeu-de-taquin} to calculate Schubert
structure constants of the cohomology of (co)minuscule flag
varieties. We do not recover this result, but we note that their
notion of standard tableau, fits into our framework as strong (or
weak) tableaux for the distributive parabolic quotients, with the
edge labels forgotten.

\section{Dual graded graphs for Kac-Moody algebras}
\label{sec:DGG}

\subsection{Dual graded graphs}
We recall Fomin's notion of dual graded graphs \cite{Fomin:dual}. A
\textit{graded graph} is a directed graph $\Gamma=(V,E,h,m)$ with
vertex set $V$ and set of directed edges $E\subset V^2$, together
with a grading function $h:V\rightarrow \Z_{\ge0}$, such that every
directed edge $(v,w)\in E$ satisfies $h(w)=h(v)+1$ and has a
multiplicity $m(v,w) \in\Z_{\ge0}$. Forgetting the edge labels $m$,
$\Gamma$ may be regarded as the Hasse diagram of a graded poset. We
shall interpret $m(v,w)$ as making $\Gamma$ into a directed
multigraph in which there are $m(v,w)$ distinct edges from $v$ to
$w$.

$\Gamma$ is \textit{locally-finite} if, for every $v\in V$, there
are finitely many $w\in V$ such that $(v,w)\in E$ and finitely many
$u\in V$ such that $(u,v)\in E$; we shall assume this condition
without further mention. For a graded graph $\Gamma=(V,E,h,m)$
define the $\Z$-linear down and up operators $D,U:\Z V\to\Z V$ on
the free abelian group $\Z V$ of formal $\Z$-linear combinations of
vertices, by
\begin{equation*}
U_\Gamma(v) = \sum_{(v,w)\in E} m(v,w) w \qquad D_\Gamma(w) =
\sum_{(v,w)\in E} m(v,w) v
\end{equation*}

A pair of graded graphs $(\Gamma,\Gamma')$ is \textit{dual} if
$\Gamma$ and $\Gamma'$ have the same vertex sets and grading
function but possibly different edge sets and edge multiplicities,
such that
\begin{align} \label{E:balance}
  D_{\Gamma'} U_\Gamma - U_\Gamma D_{\Gamma'} = r \,\mathrm{Id}
\end{align}
as $\Z$-linear operators on $\Z V$, for some fixed $r\in\Z_{>0}$. We
call $r$ the {\it differential coefficient}.  When $\Gamma =
\Gamma'$ and all the edges have multiplicity one, we obtain the
$r$-differential posets of~\cite{Sta}.

\begin{remark} \label{R:dualinfinite}
The duality property implies that $V$ is infinite.
\end{remark}

\begin{ex} \label{X:Young}
Let $\Gamma=\Y$ be Young's lattice, with $(\la,\mu)\in E$ if the
diagram of the partition $\mu$ is obtained from that of $\la$ by
adding a single cell (in which case we say that the cell is
\textit{$\la$-addable} and \textit{$\mu$-removable}), all edge
multiplicities are $1$, and $h(\la)=|\la|$ is the number of cells in
the diagram of $\la$. Then $(\Y,\Y)$ is a pair of dual graded graphs
with differential coefficient $1$.
\end{ex}

\subsection{The labeled Kac-Moody weak and strong orders}
\label{sec:KM}

In this section a new family of dual graded graphs is introduced.

Let $I$ be a set of Dynkin nodes and $A=(a_{ij})_{i,j\in I}$ be a
generalized Cartan matrix (GCM), that is, one with integer entries
which satisfies $a_{ii}=2$ for all $i\in I$, and for all $i\ne j$,
$a_{ij}\le 0$ and $a_{ij}<0$ if and only if $a_{ji}<0$. Let
$\geh=\geh(A)$ denote the Kac-Moody algebra over $\C$ associated to
$A$ \cite{Kac}, $\h\subset\geh$ the Cartan subalgebra, and $\h^*$
the dual. Let $\{\al_i\mid i\in I\} \subset \h^*$ be the simple
roots, $\{\al_i^\vee\mid i\in I\} \subset \h$ the simple coroots,
and $\{\La_i\mid i\in I\} \subset\h^*$ the fundamental weights, with
$a_{ij}=\inner{\al_i^\vee}{\al_j}$ where $\inner{\cdot}{\cdot}:
\h\times \h^* \to \C$ is the natural pairing. We assume that the
simple roots are linearly independent and the dimension of $\h$ is
chosen to be minimal. Let $W$ be the Weyl group of $\geh$: it has
generators $s_i$ for $i\in I$ and relations $s_i^2=1$ for $i\in I$
and $(s_is_j)^{m_{ij}}=1$ for $i,j\in I$ with $i\ne j$, where
$m_{ij}$ is $2,3,4,6$ or $\infty$ according as $a_{ij}a_{ji}$ is
$0,1,2,3$ or $>3$. Let $\ell:W\to \Z_{\ge0}$ be the length function
on $W$. Let $\RealRoot = W \cdot \{\alpha_i\mid i\in I\}$ be the set
of real roots and $\RealRoot^+=\RealRoot \cap \bigoplus_{i\in I}
\Z_{\ge0}\, \alpha_i$ the positive real roots. The associated coroot
$\al^\vee$ of $\al\in\RealRoot^+$ is defined by $\al^\vee = u
\alpha_i^\vee$, where $u\in W$ and $i\in I$ are such that
$\al=u\alpha_i$. For $\al\in\RealRoot$ let $s_\al = u s_i u^{-1}$
denote the reflection associated to $\al$. The \textit{strong order}
(or Bruhat order) $\le$ on $W$ is defined by the cover relations $w
\lessdot ws_\al$ whenever $\ell(ws_\al) = \ell(w) + 1$ for some $\al
\in \RealRoot^+$ and $w \in W$. The \textit{left weak order}
$(W,\preceq)$ is the subposet of $(W,\le)$ generated by the cover
relations $w \prec s_i w$ whenever $\ell(s_i w) = \ell(w)+1$ for
some $i\in I$ and $w\in W$. The \textit{left descent set} of $v$ is
defined by $\Des(v)=\{i\in I\mid s_i v \prec v \}$.

Given $\La$ in the set $P^+$ of dominant integral weights, let
$\Gs(\La)$ be the graded graph with vertex set $W$ and edges
$(v,w)\in W^2$ such that $v\lessdot w$, with multiplicity
$m_\La(v,w)=\inner{\al^\vee}{\La}$, where $\al\in\RealRoot^+$ is
such that $w=vs_\al$. Let $i\in I$ and $u\in W$ be such that
$\al=u\alpha_i$. Then
\begin{align} \label{E:stronglabel}
  m_\La(v,w) = \inner{u \alpha_i^\vee}{\La} =
  \inner{\alpha_i^\vee}{u^{-1}\La}.
\end{align}
Let $Z^+=Z^+(\geh(A)) = Z(\geh(A)) \cap \bigoplus_{i\in I} \Z_{\ge0}
\al_i^\vee$, where $Z(\geh(A))$ is the center of $\geh(A)$. If $K\in
Z^+$, writing $K=\sum_{i\in I} k_i \al_i^\vee$, the vector
$(k_i)_{i\in I}$ defines a linear dependence amongst the rows of
$A$.

Given $K\in Z^+$, let $\Gw(K)$ be the graded graph with vertex set
$W$ and edges $(v,w)\in W^2$ such that $v\prec w=s_iv$, with
multiplicity
\begin{align} \label{E:weaklabel}
  n_K(v,w)=k_i=\inner{K}{\La_i}.
\end{align}

Both $\Gs(\La)$ and $\Gw(K)$ are graded by the length function.

\begin{thm} \label{thm:main} Let $(\La,K)\in P^+\times Z^+$.
Then $(\Gs(\La),\Gw(K))$ is a pair of dual graded graphs with
differential coefficient $r=\inner{K}{\La}$.
\end{thm}
\begin{proof} Let $U=U_{\Gs(\La)}$ and $D=D_{\Gw(K)}$.
The coefficient of $u \ne v$ in $(DU-UD)v$ is given by
\begin{align*}
 \sum_{\substack{(i,\al)\in I\times \RealRoot^+ \\
 v \lessdot vs_\al \\ u=s_i vs_\al \prec vs_\al}}
k_i \inner{\al^\vee}{\La} -
 \sum_{\substack{(i,\al)\in I\times \RealRoot^+ \\
 s_i v \prec v \\ s_i v \lessdot s_ivs_\al=u}}
k_i \inner{\al^\vee}{\La}.
\end{align*}
This quantity is zero because the indexing sets of both sums
coincide, by two versions of \cite[Lemma 5.11]{H}.

For every $i\in I$ and $v\in W$, either $v\prec s_i v$ or $s_i v
\prec v$ is a covering relation. It follows that the
coefficient of $v$ in $(DU-UD)v$ is
\begin{align*}
&\sum_{i\in I\backslash\Des(v)} k_i \inner{v^{-1}\alpha_i^\vee}{\La}
  - \sum_{i\in \Des(v)} k_i \inner{(s_i v)^{-1}\alpha_i^\vee}{\La}
  \\
  &= \sum_{i\in I} k_i \inner{v^{-1}\alpha_i^\vee}{\La} = \sum_{i\in I}
  k_i
  \inner{\alpha_i^\vee}{v \La} \\
  &= \inner{K}{v \La} = \inner{v^{-1}K}{\La}=\inner{K}{\La}.
\end{align*}
We have used the $W$-invariance of $\inner{\cdot}{\cdot}$ and $K$.
\end{proof}

\begin{remark} \label{R:strongSchubert} For $i\in I$ and $v\in W$,
the multiplicity of the edge $(v,vs_\al)$ in $\Gs(\La_i)$, is the
Chevalley multiplicity, given by the coefficient of $\xi^{vs_\al}$
in the product $\xi^{s_i} \xi^v$, where $\xi^v\in H^*(G/B)$ is the
Schubert cohomology class for the flag manifold $G/B$ associated
with the Kac-Moody algebra $\geh$ \cite{KK}.
\end{remark}

In Proposition~\ref{P:homology} we will relate the multiplicities of the weak graph
$\Gw(K)$ with the homology multiplication of the affine Grassmannian, in the case that
$\geh(A)$ is of untwisted affine type.

\subsection{Tableaux and enumeration}
\label{sec:tableaux}

Let $\Gamma=(V,E,h,m)$ be a graded graph and $v, w \in V$. A
\textit{$\Gamma$-tableau} $T$ of shape $v/w$ is a directed path
$$T = (w = v_0 \stackrel{m_1}{\longrightarrow} v_1
\stackrel{m_2}{\longrightarrow} \cdots
\stackrel{m_k}{\longrightarrow} v_k = v)$$ from $w$ to $v$ in
$\Gamma$, where $m_i$ is an element taken from a set of
$m(v_{i-1},v_i)$ possible markings of the edge $(v_{i-1},v_i)$. In
the multigraph interpretation, $m_i$ indicates which of the
$m(v_{i-1},v_i)$ edges going from $v_{i-1}$ to $v_i$, is traversed
by the path.

Let $v=\outsh(T)$ and $w=\insh(T)$ the outer and inner shapes of $T$
respectively, and denote by $\Tab(\Gamma)$ the set of
$\Gamma$-tableaux, and $\Tab(\Gamma,v/w)$ the subset of those of
shape $v/w$.

\begin{ex} \label{X:stdtabs} For $\Y$ as in Example \ref{X:Young},
$\Y$-tableaux are standard Young tableaux.
\end{ex}

If $\Gamma$ has a unique minimum element $\hat{0}$, we say $T$ has
shape $v$ if $\insh(T)=\hat{0}$.

\begin{thm} \cite{Fomin:dual} Let $(\Gamma,\Gamma')$ be a pair of
dual graded graphs with differential coefficient $r$. Then
\begin{equation}
\label{E:Schenstedidentity} %
r^n\, n! = \sum_{\substack{v\in V\\ h(v)=n}} f^v_\Gamma
f^v_{\Gamma'}
\end{equation}
where $f^v_\Gamma=|\Tab(\Gamma,v)|$ and
$f^v_{\Gamma'}=|\Tab(\Gamma',v)|$.
\end{thm}

\begin{ex} \label{X:Schenstedid}
Let $\Gamma = \Gamma'=\Y$. Then by Examples \ref{X:Young} and
\ref{X:stdtabs}, equation \eqref{E:Schenstedidentity} is the
well-known identity $n!=\sum_\la f_\la^2$ where $\la$ ranges over
the partitions of $n$ and $f_\la$ is the number of standard Young
tableaux of shape $\la$.
\end{ex}

The graphs $\Gw(K)$ and $\Gs(\La)$ both have minimum element
$\idd\in W$. We call $\Gw(K)$-tableaux \textit{(standard) $K$-weak tableaux}
and $\Gs(\La)$-tableaux \textit{(standard) $\La$-strong tableaux}.

\begin{cor}
\label{cor:Fom} Let $\geh$ be a Kac-Moody algebra and $(\La,K) \in
P^+ \times Z^+$.  Then for each $n \in \Z_{\ge 0}$ we have
\begin{equation}
\label{E:identity}  r^n\, n! = \sum_{\substack{w\in W\\ \ell(w) =
n}} f_\weak^w f_\strong^w
\end{equation}
where $r=\inner{K}{\La}$, $f_\weak^w$ is the number of $K$-weak
tableaux of shape $w$ and $f_\strong^w$ is the number of
$\La$-strong tableaux of shape $w$.
\end{cor}

In section \ref{sec:LLMS} standard Young tableaux are realized as
special cases of both $K$-weak and $\La$-strong tableaux using
affine algebras of type $A^{(1)}$.

\subsection{From dual graded graphs to Schensted bijections}
\label{sec:comb} A \textit{differential bijection} for the pair of
dual graded graphs $(\Gamma,\Gamma')$ is one which exhibits the
equation \eqref{E:balance}; we shall give a precise definition
below. By \cite{Fomin:Schensted}, a differential bijection induces a
\textit{Schensted bijection} (see \eqref{E:Schenbij}) which
describes the enumerative identity \eqref{E:Schenstedidentity}.

We recall Fomin's theory in more detail. Given $r\in\Z_{\ge0}$, let
$S(r)$ be a set of cardinality $r$. Sometimes we may write $S(r)$
for a particular set of cardinality $r$. Let $(\Gamma,\Gamma')$ be a
pair of dual graded graphs with $\Gamma=(V,E,h,m)$ and
$\Gamma'=(V,E',h,m')$ and differential coefficient $r$. Given
$x,y\in V$, define $\UD_{xy}=\{(z,m,m')\in V \times \Z_{>0}^2\mid
(z,y)\in E, (z,x)\in E', m \le m(z,y), m' \le m'(z,x)\}$ and let
$\DU_{xy}=\{(w,M,M')\in V\times \Z_{>0}^2 \mid (x,w)\in E, (y,w)\in
E', M\le m(x,w), M'\le m'(y,w)\}$. $\UD_{xy}$ represents the set of
marked paths going down one step in $\Gamma'$ from $x$ to some $z\in
V$ and then up one step in $\Gamma$ from $z$ to $y$. Similarly,
$\DU_{xy}$ represents the set of marked paths going up one step in
$\Gamma$ from $x$ to some $w\in V$, and then down one step in
$\Gamma'$ from $w$ to $y$. To cancel the off-diagonal terms in
\eqref{E:balance}, for every $(x,y)\in V^2$ with $x\ne y$, there
must be a bijection
\begin{equation} \label{E:offdiagbij}
  \db_{xy}: \UD_{xy} \to \DU_{xy},
\end{equation}
and to obtain agreement of diagonal terms in \eqref{E:balance}, for
each $x\in V$ there must be a bijection
\begin{equation}\label{E:diagbij}
\db_x: S(r) \sqcup \UD_x \to \DU_x
\end{equation}
where $\DU_x=\DU_{xx}$ and $\UD_x=\UD_{xx}$. By definition, a
\textit{differential bijection} for $(\Gamma,\Gamma')$, is a
collection $\db=(\db_{xy};\db_x)$ of such bijections $\db_{xy}$ and
$\db_x$.

\begin{ex} \label{X:Youngdb} Let $\Gamma=\Y$ with dual graded graph
structure on $(\Y,\Y)$ as in Example \ref{X:Young}. For $\la\in\Y$,
since all edges have multiplicity $1$, $\DU_\la$ is in bijection
with $\la$-addable corner cells and $\UD_\la$ is in bijection with
$\la$-removable corner cells. The $\la$-addable and $\la$-removable
corner cells of $\la$ are interleaved. Let $\db_\la$ send the unique
element of the set $S(1)$, to the $\la$-addable corner cell in the
first row of $\la$, and send a $\la$-removable corner to the nearest
$\la$-addable corner with higher row index. For $\la\ne\mu$ the sets
$\DU_{\la\mu}$ and $\UD_{\la\mu}$ have the same cardinality, which
is either $0$ or $1$, so there is no choice for the definition of
$\db_{\la\mu}$. This defines a differential bijection $\db$ for
$(\Y,\Y)$.
\end{ex}

Let $P_n(r)$ be the set of $r$-colored permutations of $n$ elements.
We realize $\sigma\in P_n(r)$ as an $n\times n$ monomial matrix (one
with exactly one nonzero element in each row and in each column,
whose nonzero entries must be taken from a set $S(r)$ of cardinality
$r$ such that $0\notin S(r)$).

We assume that $\Gamma$ and $\Gamma'$ have a common minimum element
$\hat{0}$ such that $h(\hat{0})=0$.

A \textit{Schensted bijection} $\ins$ for $(\Gamma,\Gamma')$ is a
family of bijections (for all $n\in\Z_{\ge0}$)
\begin{equation}
\label{E:Schenbij}
\begin{split}
P_n(r) &\to \bigsqcup_{\substack{v\in V
\\ h(v)=n}}
\Tab(\Gamma,v) \times \Tab(\Gamma',v) \\
\sigma&\mapsto (P,Q)
\end{split}
\end{equation}
We fix a differential bijection $\db$ for $(\Gamma,\Gamma')$ and
define the induced Schensted bijection $\ins_\db$.

Given $\sigma\in P_n(r)$, we shall define a directed graph $G$ with
vertices $G_{ij}\in V$ for $0\le i,j\le n$, which is depicted
matrix-style. It shall have the property that (1) adjacent vertices
$G_{i,j-1},G_{ij}$ in a row, are either equal, or they form an edge
$(G_{i,j-1},G_{ij})\in E$ with marking $m\in
S(m(G_{i,j-1},G_{ij}))$, and (2) adjacent vertices
$G_{i-1,j},G_{ij}$ in a column, are either equal, or they form an
edge $(G_{i-1,j},G_{ij})\in E'$ with marking $m'\in
S(m'(G_{i-1,j},G_{ij}))$. Moreover, (3) $G_{i-1,j}\ne G_{ij}$ (resp.
$G_{i,j-1}\ne G_{ij}$) if and only if the unique $p$ such that
$\sigma_{pj}\ne0$ (resp. $q$ such that $\sigma_{iq}\ne0$) satisfies
$p\le i$ (resp. $q\le j$). In particular (ignoring the equalities),
for each $i$, the $i$-th row $G_{i\bullet}$ is a $\Gamma$-tableau
and for each $j$, the $j$-th column $G_{\bullet j}$ is a
$\Gamma'$-tableau. For the sake of uniform language we shall always
imagine that there is a marked edge $G_{i,j-1}\to G_{ij}$ and
$G_{i-1,j}\to G_{ij}$, but when the vertices coincide the marked
edge degenerates.

$G$ is defined inductively as follows. The north and west edges
$G_{0\bullet}$ and $G_{\bullet0}$ of $G$ are initialized to the
empty tableau: $G_{i0}=G_{0j}=\hat{0}$ for all $0\le i,j\le n$. To
define the rest of $G$, it suffices to give a \textit{local rule},
which, given the marked edges $G_{i-1,j-1}\overset{m}\to G_{i-1,j}$
and $G_{i-1,j-1}\overset{m'}\to G_{i,j-1}$, and the value
$\sigma_{ij}$, determines $G_{ij}\in V$ with markings $M\in
S(m(G_{i,j-1},G_{ij}))$ and $M'\in S(m'(G_{i-1,j},G_{ij}))$.

This is depicted below.  Use $z,y,x,w$ to denote
$G_{i-1,j-1},G_{i-1,j},G_{i,j-1},G_{ij}$ for convenience and write
$c=\sigma_{ij}$. In later examples we shall indicate $\sigma_{ij}=1$
by the symbol $\otimes$ and $\sigma_{ij}=0$ by a blank.
\begin{equation}\label{E:fominsquare}
\xymatrix@=5pt{%
{G_{i-1,j-1}} \ar[rr]^m \ar[dd]_{m'} && {G_{i-1,j}} \ar@{-->}[dd]^{M'} \\
& {\sigma_{ij}}
\\
{G_{i,j-1}} \ar@{-->}[rr]_M && {G_{ij}} %
} %
\qquad\quad
\xymatrix@=5pt{%
{z} \ar[rr]^m \ar[dd]_{m'} && {y} \ar@{-->}[dd]^{M'} \\ & {c}
\\
{x} \ar@{-->}[rr]_M && {w} %
}
\end{equation}
The local rule is defined using $\db$.
\begin{enumerate}
\item
If $z=x=y$:
\begin{enumerate}
\item
If $c=0$, set $w=z$.
\item
If $c\ne 0$, let $\db_x(c)=(w,M,M')$.
\end{enumerate}
\item
If $z\ne x=y$ then let $\db_x(z,m,m')=(w,M,M')$.
\item
If $z=y$ and $z\ne x$ then let $w=x$ and $M'=m'$.
\item
If $z=x$ and $z\ne y$ then let $w=y$ and $M=m$.
\item
If $z,x,y$ are all distinct, then let $(w,M,M')=\db_{xy}(z,m,m')$.
\end{enumerate}
This uniquely determines $G$ \cite{Fomin:Schensted}. Its south edge
$G_{n\bullet}$ is a $\Gamma$-tableau $P$ and its east edge
$G_{\bullet n}$ is a $\Gamma'$-tableau $Q$, both of a common shape
$v=G_{nn}\in V$ with $h(v)=n$. This well-defines a map $\ins_\db$ as in
\eqref{E:Schenbij}.

For the inverse of $\ins_\db$, let $v\in V$ be such that $h(v)=n$,
and let $(P,Q)\in \Tab(\Gamma,v)\times \Tab(\Gamma',v)$. To recover
$\sigma\in P_n(r)$, we initialize the south and east edges of $G$ to
$P$ and $Q$ respectively. Then for each $i,j$ and two by two
subgraph as above, we apply the inverse of the above local rule.
Given labeled edges $x\overset{M}\to w$ and $y\overset{M'}\to w$, it
determines $z\in V$ and marked edges $z\overset{m}\to y$ and
$z\overset{m'}\to x$ and a value $c\in \{0\}\sqcup S(r)$, such that
$c\ne 0$ if and only if $z=x=y\ne w$ and $\db_z^{-1}(w,M,M')=c$. The
inverse local rule is defined as follows.
\begin{enumerate}
\item
If $x=y$:
\begin{enumerate}
\item
If $w=x$, let $z=x$.
\item
If $w\ne x$:
\begin{enumerate}
\item If $c:=\db_x^{-1}(w,M,M')\in S(r)$:
let $z=x$.
\item Otherwise $\db_x^{-1}(w,M,M')=(z,m,m')$.
\end{enumerate}
\end{enumerate}
\item If $w=x\ne y$, let $z=y$ and $m'=M'$.
\item If $w=y\ne x$, let $z=x$ and $m=M$.
\item If $x,y,w$ are all distinct,
let $(z,m,m')=\db_{xy}^{-1}(w,M,M')$.
\end{enumerate}
In all cases but (1)(b)(i) let $c=0$. Using the inverse local rule
the rest of $G$ is defined \cite{Fomin:Schensted} and one obtains a
well-defined element $\sigma\in P_n(r)$.

\begin{thm} \cite{Fomin:Schensted} Let $(\Gamma,\Gamma')$ be a dual
graded graph with differential coefficient $r$. Then for any
differential bijection $\db$ for $(\Gamma,\Gamma')$, the above
construction defines a Schensted bijection $\ins_{\db}$ of the form
\eqref{E:Schenbij}.
\end{thm}

We say that $\ins_{\db}$ is the \textit{Schensted bijection} induced
by the differential bijection $\db$.

\begin{ex} \label{X:YoungSchen}
The differential bijection $\db$ of Example \ref{X:Youngdb} induces
Schensted's row insertion bijection \cite{Schensted}.
\end{ex}

\begin{remark} \label{R:KMoffdiag}
For the Kac-Moody dual graded graphs $(\Gw(K),\Gs(\La))$ there is a
natural choice for the off-diagonal part of the differential
bijection. If $v\ne w$ with $v,w\in W$ then $\db_{vw}$ is
essentially obtained from \cite[Lemma 5.11]{H}, just as in the proof
of Theorem \ref{thm:main}. For $i\in \Des(v)$, the marked
``down-then-up'' path $v \stackrel{m'}{\longrightarrow} s_iv
\stackrel{m}{\longrightarrow} s_ivs_\alpha = w$ maps to the marked
``up-then-down'' path $v \stackrel{m}{\longrightarrow} vs_\alpha
\stackrel{m'}{\longrightarrow} s_ivs_\alpha = w$. Here $m$ and $m'$
denote edge markings, which in either case are selected from sets of
size
$m_{\La}(s_iv,s_ivs_\alpha)=\inner{\al^\vee}{\La}=m_{\La}(v,vs_\alpha)$
and $n_K(\La_i)$ respectively.

Currently we are not aware of a general rule for $\db_v$ which
exhibits the coefficient of $v$ in $(DU-UD)v$ as $\inner{K}{\La}$.
In section \ref{sec:LLMS} we shall give a special case where the
bijection $\db_v$ has been constructed explicitly.
\end{remark}

\subsection{Automorphisms and mixed insertion}
\label{sec:mixed} This section is a natural synthesis of the ideas
of Fomin \cite{Fomin:dual} and Haiman \cite{Ha} which does not seem
to have been written down before.  We believe this construction is
particularly interesting for Kac-Moody dual graded graphs (see also
Section \ref{sec:fold}).

Let $(\Gamma,\Gamma')$ be a pair of dual graded graphs with
$\Gamma=(V,E,m,h)$ and $\Gamma=(V,E',m',h)$. Say that a permutation
$\tau:V\to V$ is an automorphism of $(\Gamma,\Gamma')$ if (1)
$h\circ \tau = h$, (2) $(x,y)\in E$ if and only if
$(\tau(x),\tau(y))\in E$, and in this case,
$m(x,y)=m(\tau(x),\tau(y))$, and (3) $(x,y)\in E'$ if and only if
$(\tau(x),\tau(y))\in E'$, and in this case,
$m'(x,y)=m'(\tau(x),\tau(y))$.

Given a differential bijection $\db$ for $(\Gamma,\Gamma')$, we
define its twist $\db^\tau$ by $\tau$ as follows. For every $x,y\in
V$ there are natural bijections $\tau:\UD_{xy}\to
\UD_{\tau(x)\tau(y)}$ given by $(z,m,m')\mapsto (\tau(z),m,m')$ and
$\tau:\DU_{xy}\to\DU_{\tau(x)\tau(y)}$ given by $(w,M,M')\mapsto
(\tau(w),M,M')$. Let $\tau: S(r)\to S(r)$ be the identity. Then
define $\db^\tau_{xy}= \tau^{-1} \circ \db_{\tau(x)\tau(y)} \circ
\tau$. It is easy to verify that $\db^\tau$ is also a differential
bijection for $(\Gamma,\Gamma')$.

\begin{ex} \label{X:Schenauto}
Let $\Gamma=\Gamma'=\Y$ and $\tr:\Y\to\Y$ the automorphism of
$(\Gamma,\Gamma)$ that transposes partition diagrams. Let $\db$ be
the differential bijection in Example \ref{X:Youngdb}. Then
$\ins_{\db^\tr}$ is Schensted's column insertion bijection
\cite{Schensted}.
\end{ex}

For the sequel we assume that $\tau$ has finite order $\kappa$. A
\textit{$\tau$-mixed} $\Gamma$-tableau is a $\Gamma$-tableau whose
edges have an auxiliary marking parameter $p\in
S(\kappa)=\{0,1,\dotsc,\kappa-1\}$. Let $\Tab_\tau(\Gamma)$ be the
set of $\tau$-mixed $\Gamma$-tableaux. Suppose $\tau'$ is an
automorphism of $(\Gamma,\Gamma')$ of order $\kappa'$. Let
$(\Gamma,\Gamma';\tau,\tau')$ denote the pair of dual graded graphs
given by $\Gamma$ and $\Gamma'$ except that $\Gamma$-edges (resp.
$\Gamma'$-edges) are labeled by $(m,p)$ with $p\in S(\kappa)$ (resp.
$(m',p')$ with $p'\in S(\kappa')$) and $m$ (resp. $m'$) is a usual
edge label for $\Gamma$ (resp. $\Gamma'$). This multiplies the
number of markings for each edge of $\Gamma$ (resp. $\Gamma'$) by
$\kappa$ (resp. $\kappa'$). The differential coefficient of
$(\Gamma,\Gamma';\tau,\tau')$ is $r\kappa\kappa'$ where $r$ is the differential
coefficient of $(\Gamma,\Gamma')$.

Let $\db$ be a differential bijection for $(\Gamma,\Gamma')$. Then
there is an obvious differential bijection (also denoted $\db$) for
$(\Gamma,\Gamma';\tau,\tau')$, defined by a trivial scaling by
$\kappa$ in $\Gamma$ and by $\kappa'$ in $\Gamma'$.

We shall define another bijection
\begin{align} \label{E:mixedbij}
  P_n(\kappa\kappa' r) &\to \bigsqcup _{\substack{x\in V\\ h(x)=n}}
  \Tab_\tau(\Gamma,x) \times \Tab_{\tau'}(\Gamma',x).
\end{align}
called \textit{$(\tau,\tau')$-mixed insertion} by modifying the
process in which we construct the matrix $G_{ij}$ from the colored
permutation $\sigma$. Instead of using the same differential
bijection $\db$ to compute each $G_{ij}$, we use twists of $\db$ by
automorphisms that depend on $(i,j)$.

Let $\sigma\in P_n(\kappa\kappa' r)$. We regard $\sigma$ as a
monomial matrix in which each nonzero entry has three labels
$(c,p,p')\in S(r)\times S(\kappa)\times S(\kappa')$, where
$S(\kappa)=\{0,1,\dotsc,\kappa-1\}$ and
$S(\kappa')=\{0,1,\dotsc,\kappa'-1\}$. Each horizontal (resp.
vertical) edge is marked by a pair $(m,p)$ (resp. $(m',p')$) where
$m$ (resp. $m'$) is the usual marking and $p\in S(\kappa)$ (resp.
$p'\in S(\kappa')$). Let
$(z,y,x,w)=(G_{i-1,j-1},G_{i-1,j},G_{i,j-1},G_{ij})$ and suppose
that $z \overset{(m,p)}\to y$ and $z\overset{(m',p')}\to x$ are
given, where it is understood that if $z=y$ (resp $z=x$) then
$(m,p)$ (resp. $(m',p')$) need not be specified. Then $G_{ij}$ is
determined as before, except that instead of using $\db_{xy}$ we use
the twist $\db_{xy}^{\tau^k (\tau')^{k'}}$ where $k$ and $k'$ are as
follows.  Let $(c,p,p')$ be the nonzero entry of $\sigma$ in the
$i$-th row, say, $\sigma_{il}$.  We set $k' = p'$.  Separately, let $(c,p,p')$ be the nonzero entry of
$\sigma$ in the $j$-th column, say, $\sigma_{qj}$.  We set $k = p$.
Note that in the case $q > i$  (resp. $l > j$) the bijection $\db_{xy}$ is not used in the local
rule so the value of $k$ (resp. $k'$) does not affect the algorithm.


In other words,
if $\sigma_{ij}=(c,p,p')$ is a nonzero entry, then $(\tau')^{p'}$
acts everywhere to the right in the $i$-th row and all vertical
edges to the right (those of the form $G_{i-1,l}\to G_{il}$ for
$l\ge j$) are given the auxiliary marking $p'$, and $\tau^p$ acts
everywhere below in the $j$-th column, and all horizontal edges
below (those of the form $G_{l,j-1}\to G_{lj}$ for $l\ge i$) are
given the auxiliary marking $p$. The output is the pair $(P,Q)\in
\Tab_\tau(\Gamma,v)\times \Tab_{\tau'}(\Gamma',v)$ where $v=G_{nn}$
and $P$ and $Q$ are obtained from the south and east edges of $G$
respectively.

\begin{prop} \label{P:mixed} $(\tau,\tau')$-mixed insertion gives a
well-defined bijection \eqref{E:mixedbij}.
\end{prop}

\begin{ex} \label{X:mixed}
In the context of Example \ref{X:Schenauto}, $(\tau,\tau')$-mixed
insertion specializes to the following kinds of insertion
algorithms.
\begin{enumerate}
\item $(\tau,\tau')=(\idd,\idd)$: Schensted row insertion
\cite{Schensted}.
\item $(\tau,\tau')=(\idd,\tr)$: left-right insertion \cite{Ha}.
\item $(\tau,\tau')=(\tr,\idd)$: mixed insertion \cite{Ha}.
\item $(\tau,\tau')=(\tr,\tr)$: doubly-mixed insertion \cite{Ha}.
\end{enumerate}
\end{ex}

\subsection{Restriction to parabolics}
Let $J \subset I$.  We say that a weight $\La$ is supported on $J$
if $\La = \sum_{j \in J} a_j \La_j$.  The Kac-Moody dual graded
graphs $(\Gs(\La),\Gw(K))$ are compatible with restriction to
parabolics.  Let $W_J\subset W$ be the parabolic subgroup generated
by $\{s_j\mid j\in J\}$ and let $W^J$ be the set of minimal length
coset representatives in $W/W_J$.  Note that $W^J$ inherits weak and
strong orders from $W$ by restriction.

\begin{prop}\label{P:parabolic}
Fix $J \subset I$.  If $\Lambda$ is supported on $I\backslash J$
then the restriction of $(\Gs(\La),\Gw(K))$ to $W^J$ is a pair of
dual graded graphs with differential coefficient $\inner{K}{\La}$.
\end{prop}
\begin{proof}
Suppose $v \prec w$ is a weak cover. If $w \in W^J$ then $v \in W^J$
since $W^J \subset W$ is a lower order ideal for $\preceq$.

Suppose $w \lessdot v$ and $w \in W^J$ and $v \notin W^J$.  Since
$v$ has a reduced expression ending in $s_j$ for some $j \in J$ and
$w$ is obtained from this reduced expression by omitting a simple
generator, we conclude that $v = ws_j$.  But $\Lambda$ is supported
on $I\backslash J$, so $\inner{\alpha_j^\vee}{\La} = 0$.

Combining these two facts we see that the proof of
Theorem~\ref{thm:main} restricts to $W^J$.
\end{proof}

We shall use the following notation for maximal parabolic subgroups
of $W$. For $i\in I$ we shall write $W^i$ for $W^J$ where
$J=I\setminus\{i\}$. We denote by $(\Gs(\La),\Gw(K))^i =
(\Gs^i(\La),\Gw^i(K))$, the dual graded graph given by restriction
of $(\Gs(\La),\Gw(K))$ to $W^i$.

\subsection{The affine case}
\label{sec:affine} If the GCM $A$ is of finite type, then
$Z^+=\{0\}$ and all of the edges of $\Gw(K)$ are labeled $0$.

In this section let $A$ be of untwisted affine type. Let $0\in I$ be
the distinguished Kac $0$ node and $J=I\setminus\{0\}$. Then $W$ is
the affine Weyl group, $W_J=W_{\rm fin}$ is the finite Weyl group,
and we write $W^0=W^J$. By Proposition \ref{P:parabolic} the
restriction of the Kac-Moody dual graded graph to $W^0$, is a dual
graded graph. In this case the weak graph $\Gw(K)$ has an
interpretation involving the Schubert calculus of the homology of
the affine Grassmannian, and the duality is a combinatorial
expression of the pairing between the homology and cohomology of the
affine Grassmannian.

For affine algebras, $Z^+=\Z_{\ge0}K$ where $K = K_{\rm can} =
\sum_{i \in I}k_i \alpha_i^\vee$ is the canonical central element;
the vector $(k_i)_{i\in I}$ is the unique linear dependence of the
rows of $A$ given by positive relatively prime integers \cite{Kac}.
In this case, since the labels of $\Gw(K)$ are linear in $K$,
without loss of generality we shall only work with $K_{\rm can}$ and
define $\Gw := \Gw(K_{\rm can})$.
%
The edge labels of $\Gw$ are related to the homology multiplication
in affine Grassmannians, as follows.

Let $\geh = \geh(A)$ be an untwisted affine algebra.  Let
$\rm{Gr}=\rm{Gr}_G$ denote the affine Grassmannian of the simple Lie
group $G$ whose Lie algebra $\geh_{\rm fin}$ is the canonical simple
Lie subalgebra of the affine algebra $\geh$. For $w \in W^0$ we let
$\xi_w \in H_*({\rm Gr})$ denote the corresponding {\it homology}
Schubert class. Recall the constants $n(w,v)$ from
\eqref{E:weaklabel}.

\begin{prop}\label{P:homology}
Let $\xi_0 = \xi_{s_0}$ be the Schubert class indexed by the unique
simple generator $s_0 \notin W_{\rm fin}$. Then for every $w\in
W^0$, we have in $H_*({\rm Gr})$ the identity
$$
\xi_0 \; \xi_w = \sum_v n(w,v) \,\xi_v
$$
where $v\in W^0$ runs over the weak covers $w\prec v$ of $w$.
\end{prop}
\begin{proof}
We rely on the results of~\cite{Lam} which in turn are based on
unpublished work of Peterson.  Let $S = {\rm Sym}((\h^*_\Z)_\fin) =
H^{T_\fin}({\rm pt})$ denote the symmetric algebra in the weights of
the $\geh_{\rm fin}$ and $\phi_0: S \to \Z$ denote the evaluation at
0.  Let $\A_0$ denote the affine nilCoxeter algebra corresponding to
$W$.  As a free $\Z$-module $\A_0$ is spanned by elements $\{A_w
\mid w \in W\}$ with multiplication given by
$$
A_w \, A_v =\begin{cases} A_{wv} & \mbox{if $\ell(w) + \ell(v) =
\ell(wv)$} \\
0 & \mbox{otherwise.} \end{cases}
$$
The affine nilHecke algebra $\A$ is the $\Z$-algebra generated by
$\A_0$ and $S$ with the addition relation (\cite[Lemma 3.1]{Lam}):
\begin{equation}\label{eq:commute}
A_w \, \lambda = (w \cdot \lambda)A_w + \sum_{w\,r_\alpha  \lessdot
w}\inner{\lambda}{\alpha^\vee} A_{w\,r_\alpha},
\end{equation}
where $\alpha$ is always taken to be a positive root of $W$.

Now let $$ \B = \{a \in \A_0 \mid \phi_0(as) =\phi_0(s)a \;
\text{for any} \; s \in S \} \subset \A_0
$$
denote the affine Fomin-Stanley subalgebra, where $\phi_0:\A \to
\A_0$ is given by $\phi_0(\sum_{w} a_w \,A_w) = \sum_w \phi_0(a_w)
A_w$. Let $j_0: H_*(\Gr_G) \to \B$ denote the ring isomorphism
(\cite[Theorem 5.5]{Lam}) from the homology of $\Gr_G$ to the affine
Fomin-Stanley algebra $\B$.  We first show that $j_0(\xi_0) =
\sum_{i \in I} k_i A_i$, where $A_i$ are the generators of the
nilCoxeter algebra and $K_{\rm can} = \sum_i k_i \alpha_i^\vee$.

By \cite[Proposition 5.4]{Lam}, the element $j_0(\xi_0)$ is
characterized by having unique Grassmannian term $A_0$, and the
property that it lies in $\B$.  Since $k_0 = 1$, the unique
Grassmannian term property is immediate. Using \eqref{eq:commute},
we calculate that
$$
\phi_0(a \alpha_j) = \sum_{i \in I}
k_i\inner{\alpha_i^\vee}{\alpha_j} = \inner{K_{\rm can}}{\alpha_j} =
0.
$$
Since the $\alpha_j$ span $(\h^*_\Z)_\fin$ over $\Q$, we deduce that
$\phi_0(as) = \phi_0(s) a$ for any $s \in S$, and thus $j_0(\xi_0) =
\sum_{i \in I} k_i A_i$.\footnote{This observation was first pointed
out to us by Alex Postnikov.}

By \cite[Lemma 4.3, Theorem 5.5]{Lam}, we thus have
\begin{align*}
\xi_0 \; \xi_w &= j(\xi_0) \cdot \xi_w =(\sum_{i \in I} k_i A_i)
\cdot  \xi_w  = \sum_{w \prec s_i\, w } k_i \, \xi_{s_i\,w}.
\end{align*}
where we have used the action of $\A_0$ on $H_*(\Gr_G)$ given by
(see \cite[(3.2)]{Lam})
$$
A_i \cdot \xi_w = \begin{cases} \xi_{s_i\, w} & \mbox{ if $s_iw >
w$} \\ 0 & \mbox{otherwise.} \end{cases}
$$

Recalling the definition $n(w,s_iw) = k_i$ of the weak graph $\Gw$
from \eqref{E:weaklabel} this completes the proof.
\end{proof}

\section{Affine type A and LLMS insertion}
\label{sec:LLMS}

For this section let $\geh(A)$ be the affine algebra of type
$A_{n-1}^{(1)}$. In this case the combinatorics of the pair of dual
graded graphs $(\Gs(\La_i),\Gw)$ was studied extensively
in~\cite{LLMS} and was one of the main motivations of the current
work. The affine insertion algorithm of \cite{LLMS} (which we shall
call LLMS insertion) furnishes an explicit differential bijection
for $(\Gs(\La_i),\Gw)$. LLMS insertion involves nontrivial
extensions of the notion of tableaux to semistandard weak and strong
tableaux, and proves Pieri rules (formulae for certain Schubert
structure constants) in the homology $H_*({\rm Gr})$ and cohomology
$H^*({\rm Gr})$ of the affine Grassmannian of $SL(n,\C)$
\cite{LLMS}.

For type $A_{n-1}^{(1)}$ the coefficients of the canonical central
element $K$ are all $1$. Therefore the weak graph $\Gw$ has all edge
multiplicities equal to $1$. Using the rotational symmetry of the
Dynkin diagram $A_{n-1}^{(1)}$, we may assume that $\Lambda =
\Lambda_0$ and for brevity we write $\Gs$ for $\Gs(\La_0)$.

Let $I=\{0,1,\dotsc,n-1\}$ and let the Cartan matrix be defined by
$a_{i,i+1}=a_{i+1,i}=-1$ for all $i$, with indices taken modulo $n$,
$a_{ii}=2$ for all $i\in I$, and $a_{ij}=0$ otherwise. As in section
\ref{sec:KM} the Weyl group is defined by $m_{i,i+1}=3$ and
$m_{ij}=2$ for $|i-j|\ge 2$.

\subsection{Affine permutations}
We use the following explicit realization of the affine symmetric
group $W=\tS$. A bijection $w: \Z \to \Z$ is an affine permutation
with period $n$ if $w(i + n) = w(i)$ for each $i \in \Z$ and
$\sum_{i = 1}^{n} (w(i)-i)=0$. The set of affine permutations with
period $n$ form a group isomorphic to $\tS$, with multiplication
given by function composition. The reflections $t_{ij}$ in $\tS$ are
indexed by a pair of integers $(i,j)$ satisfying $i<j$ and $i \neq j
\mod n$. Suppose $v \lessdot vt_{ij} = w$ is a cover in $\tS$. Then
the edge $(v,w)$ in $\Gs$ has multiplicity equal to $\#\{k\in\Z\mid
\text{$v(i) \le k < v(j)$ and $k=0\mod n$ }\}$ \cite{LLMS}.

\subsection{Action of $\tS$ on partitions}
Given a partition $\la$, one may associate a bi-infinite binary word
$p(\la)=p=\dotsm p_{-1} p_0 p_1 \dotsm $ called its \textit{edge
sequence}. The edge sequence $p(\la)$ traces the border of the
(French) diagram of $\la$, going from northwest to southeast, such
that every letter $0$ (resp. $1$) represents a south (resp. east)
step, and some cell in the $i$-th diagonal is touched by the steps
$p_{i-1}$ and $p_i$. Here the cell $(i,j)$ lies in row $i$ (where
row indices increase from south to north), column $j$ (where column
indices increase from west to east), and diagonal $j-i$.

The affine symmetric group $\tS$ acts on partitions, since elements
of $\tS$ are certain permutations $\Z\to \Z$ and partitions can be
identified with their edge sequences, which are certain functions
$\Z\to \{0,1\}$. Then for $i\in\Z/n\Z$, $s_i\la$ is obtained by
removing from $\la$ every $\la$-removable cell of residue $i$, and
adding to $\la$ every $\la$-addable cell of residue $i$. Here the
residue of a box $(i,j)$ is the diagonal index $j-i$ taken modulo
$n$.

\subsection{Cores and affine Grassmannian permutations}
Using the language of cores, we shall describe the combinatorics of
the dual graded graph $(\Gs,\Gw)^0$ afforded by Proposition
\ref{P:parabolic}.

An \textit{$n$-ribbon} is a skew partition diagram $\la/\mu$ (the
difference of the diagrams of the partitions $\la$ and $\mu$)
consisting of $n$ rookwise connected cells, all with distinct
residues. We say that this ribbon is \textit{$\la$-removable} and
\textit{$\mu$-addable}. An \textit{$n$-core} is a partition that
admits no removable $n$-ribbon.  Since the removal of an $n$-ribbon
is the same thing as exchanging bits $p_i=0$ and $p_{i+n}=1$ in the
edge sequence for some $i$, it follows that $\la$ is a core if and
only if for every $i$, the sequence $p^{(i)}(\la) := \dotsm p_{i-2n}
p_{i-n} p_i p_{i+n} p_{i+2n}\dotsm$ consisting of the subsequence of
bits indexed by $i$ mod $n$, has the form $\dotsm 1111100000\dotsm$.
We denote the set of $n$-cores by $\Core$.

\begin{prop}[\cite{LLMS,MM}]
\label{prop:LLMS} The map $w \mapsto w \cdot \vn$ is a bijection $c:
\tS^0 \to \Core$.  Moreover, for $v, w \in \tS^0$, we have $v \le w$
if and only if $c(v) \subseteq c(w)$, and if $v \lessdot w$ then
$c(w)/c(v)$ is a disjoint union of translates of some ribbon $R$,
and the number of components of $c(w)/c(v)$ is equal to the
multiplicity $m(v,w)$ in $\Gs$.
\end{prop}

We say that $\mu \in \Core$ covers $\lambda \in \Core$ if
$c^{-1}(\mu) \gtrdot c^{-1}(\lambda)$.  Thus a standard strong
tableau in $\Gs$ is a sequence $\lambda = \lambda^0 \subset
\lambda^1 \subset \cdots \subset \lambda^l = \mu$ such that
$\lambda^{i}$ covers $\lambda^{i-1}$ and one of the components of
$\lambda^{i}/\lambda^{i-1}$ has been marked.

It is easy to show that a core cannot have both an addable and a
removable cell of the same residue. Thus, in the special case that
$v \lessdot s_iv = w$ for some $i\in I$, $c(w)$ is obtained from
$c(v)$ by adding all $c(v)$-addable cells of residue $i$, and the
ribbon $R$ of Proposition~\ref{prop:LLMS} must be a single box. In
this case we say that $c(v) \subset c(w)$ is a weak cover.

\subsection{LLMS insertion}
In \cite{LLMS}, for the affine symmetric group $\tS$, semistandard
analogues of weak and strong tableaux were defined (for all of
$\tS$, not just $\tS^0$), and an RSK correspondence was given
between a certain set of biwords or matrix words, and pairs of
tableaux, one semistandard weak and the other semistandard strong.
Let us consider the following restriction of this bijection. We
first restrict to ``standard" tableau pairs, that is, the case in
which the tableaux are weak and strong tableaux as defined in
section \ref{sec:tableaux}. Next we take the parabolic restriction
from $\tS$ to $\tS^0$. Let us denote the restricted bijection by
$I_{\rm LLMS}$.

Let $\db$ be the differential bijection for $(\Gs,\Gw)^0$ such that
$\ins_\db = I_{\rm LLMS}$. We describe it explicitly.

For $u,v\in \tS^0$ with $u\ne v$, the off-diagonal part $\db_{uv}$
of $\db$ coincides with the natural definition given in Remark
\ref{R:KMoffdiag}. The diagonal part $\db_v$ for $v\in \tS^0$ is
specified as follows. Let $\lambda = c(v)$ be the $n$-core
corresponding to $v$. If $\lambda \subset \mu$ is a weak cover then
$\mu/\lambda$ consists of all the $\la$-addable cells of $\lambda$
which have a fixed residue. As $\mu$ varies over all the weak covers
of $\lambda$ we obtain all the $\la$-addable cells in this way. Thus
there is a natural identification of the set $\DU_v$ with the set of
$\la$-addable cells. Similarly $\UD_v$ may be identified with the
set of $\la$-removable cells. This given, we may use the
differential bijection denoted $\db_\la$ in Example \ref{X:Youngdb}
for Young's lattice. This defines a differential bijection $\db$ for
$(\Gs,\Gw)^0$.

\begin{ex} \label{X:LLMS}
Figure~\ref{F:growthLLMS} shows the calculation of $I_{\rm LLMS}$ of
the permutation $\sigma=412635$ (written here in one-line notation;
it corresponds to the permutation matrix with ones located positions
$(i,\sigma(i))$ for $1\le i\le 6$) for $\tilde{S}_3$. The symbols
$\otimes$ encode $\sigma$ as described above equation
\eqref{E:fominsquare}. Each arrow indicates a marked strong cover;
the subscript $2$ indicates that the marked component is the second
from the southeast, and no subscript means the marked component is
the southeastmost. Stars in the $P$ tableau indicate the marked
components.
\Yboxdim{4pt}%
\begin{figure}
\begin{align*}
\xymatrix@=4pt{%
{.} && {.} && {.} && {.} && {.} && {.} && {.}  \\
&&&&&&&{\otimes} \\
{.} && {.}&& {.} && {.} \ar[rr]&& {\yng(1)} &&
{\yng(1)}&& {\yng(1)}  \\
&{\otimes} \\
{.} \ar[rr]&& {\yng(1)}  && {\yng(1)} && {\yng(1)} \ar[rr]&&{\yng(1,1)} &&{\yng(1,1)}
 && {\yng(1,1)} \\
&&&{\otimes} \\
{.} \ar[rr]&& {\yng(1)}\ar[rr] && {\yng(2)} && {\yng(2)} \ar[rr]&&{\yng(1,1,2)} &&{\yng(1,1,2)}
&& {\yng(1,1,2)} \\
&&&&&&&&&&&{\otimes} \\
{.} \ar[rr]&& {\yng(1)}\ar[rr] && {\yng(2)} && {\yng(2)} \ar[rr]&&{\yng(1,1,2)} &&{\yng(1,1,2)}
\ar[rr]&& {\yng(1,1,3)}\\
&&&&&{\otimes} \\
{.} \ar[rr]&& {\yng(1)}\ar[rr] && {\yng(2)} \ar[rr] && {\yng(1,3)}
\ar[rr] &&{\yng(1,1,3)} &&{\yng(1,1,3)}
\ar[rr]_{2}&& {\yng(1,1,2,4)} \\
&&&&&&&&&{\otimes} \\
{.} \ar[rr]&& {\yng(1)}\ar[rr] && {\yng(2)} \ar[rr]&& {\yng(1,3)}
\ar[rr] &&{\yng(1,1,3)} \ar[rr]&&{\yng(1,1,2,4)} \ar[rr] &&
{\yng(1,1,2,2,4)}
}%
\end{align*}
\newcommand{\prb}{2'}
\newcommand{\pre}{5'}
\newcommand{\prf}{6'}
\Yboxdim{10pt}
\begin{align*}
P &= \tableau{ 6\\ 5 \\ 4^*&6^* \\ 3&5 \\ 1^* &2^*&3^*&5^*} & Q &=
\tableau{ 6\\5\\3&6\\2&5\\1&3&4&5}
\end{align*}
\caption{Growth diagram for LLMS insertion of $412635$ for
$\tilde{S}_3$.} \label{F:growthLLMS}
\end{figure}
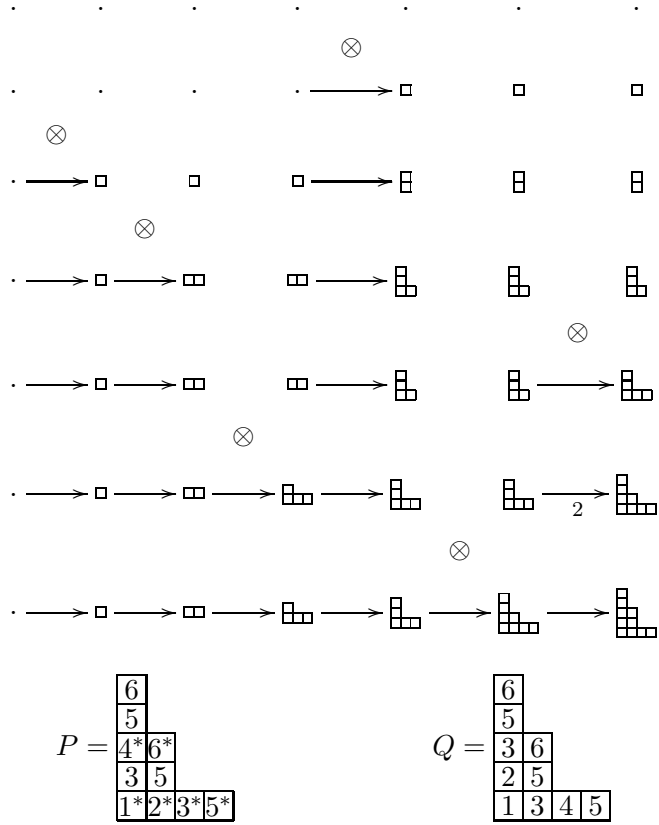
\end{ex}

\begin{remark} \label{R:LLMSlargerank}
As $n$ goes to infinity, $(\Gs,\Gw)^0$ converges to the dual graded
graph $(\Y,\Y)$ of Example \ref{X:Young} and LLMS insertion
converges to Schensted row insertion \cite{LLMS} because the
respective differential bijections coincide in the limit.
\end{remark}

\section{Folding}
\label{sec:fold}

An automorphism of the GCM $B=(b_{ij}\mid i,j\in J)$ is a
permutation $\pi$ of $J$ such that $b_{\pi(i)\pi(j)}=b_{ij}$ for all
$i,j \in J$. The automorphism $\pi$ is \textit{admissible} if
$b_{ij}=0$ for all $i$ and $j$ in the same $\pi$-orbit.

A GCM $A=(a_{ij}\mid i,j\in I)$ is \textit{symmetric} if it is a
symmetric matrix. It is \textit{symmetrizable} if there are positive
integers $o_i$ for $i\in I$, such that $DA$ is symmetric, where $D$
is the diagonal matrix with diagonal entries $o_i$.

Lusztig \cite{L} showed that every symmetrizable GCM $A$ can be
constructed from a symmetric GCM $B$ that is equipped with an
admissible automorphism $\pi$. We call this construction ``folding".

For $A$ and $B$ related in this manner, we show that the structure
of every dual graded graph of the form $(\Gs^A(\La_{i'}),\Gs^A(K))$
for $\geh(A)$, is encoded by some dual graded graph for $\geh(B)$.
Thus the combinatorics for $\geh(A)$ is reduced to that of
$\geh(B)$. In particular, for any affine algebra $\geh(A)$ there is
a simply-laced affine algebra $\geh(B)$ related by folding, so that
all affine Schensted bijections can be realized using only the
simply-laced affine algebras.

\subsection{Folding data}
Let $B=(b_{ij}\mid i,j\in J)$ be a GCM and $\pi$ an admissible
automorphism of $B$. Let $I$ be a set which indexes the $\pi$-orbits
of $J$; we write $O_i\subset J$ for the $\pi$-orbit indexed by $i\in
I$. Let $o_i=|O_i|$. It is easy to show that the matrix
$A=(a_{i'i}\mid i',i\in I)$ defined by
\begin{align} \label{E:asfold}
  a_{i'i} = \dfrac{o_{i'}}{o_i} \sum_{j\in O_i}
  b_{j'j}\qquad\text{for any $j'\in O_{i'}$},
\end{align}
is a well-defined GCM. We say that $A$ is obtained from $(B,\pi)$ by
folding.

\begin{prop} \cite[Prop. 14.1.2]{L}
Given any symmetrizable GCM $A$, there is a symmetric GCM $B$ with
admissible automorphism $\pi$, such that $A$ is obtained from
$(B,\pi)$ by folding. In particular, if $A$ is of affine type then
$B$ can be taken to be of simply-laced affine type.
\end{prop}

For the Kac-Moody algebras $\geh(A)$ and $\geh(B)$, we denote their
weight lattices by $P_A$ and $P_B$, their coweight lattices by
$P^\vee_A$ and $P^\vee_B$ and their coroot lattices by $Q^\vee_A$
and $Q^\vee_B$. For simplicity, we let $\Lambda_i, \alpha_i,
\alpha_i^\vee$ be the fundamental weights, simple roots, and simple
coroots for $\geh(A)$ and write $\omega_i, \beta_i, \beta_i^\vee$
for the corresponding data for $\geh(B)$.

Let $P'_A = P_A/(\bigoplus_{i\in I} \Z \alpha_i^\vee)^0$ be the
weight lattice of $A$ modulo the annihilator of the coroots
$\{\alpha_i^\vee\}$. Similarly define $P'_B$.  Note that
\begin{align}
\label{E:root2weightA} %
\alpha_i &= \sum_{i'\in I} a_{i'i}\,\Lambda_{i'} \\
\label{E:root2weightB} %
\beta_j &=\sum_{j'\in J} b_{j'j}\,\omega_{j'}.
\end{align}
where $\alpha_i,\Lambda_i$ and $\beta_j,\omega_j$ also denote their
respective images inside $P'_A$ and $P'_B$. Define
$$\kappa=\text{lcm}_{i\in I}(o_i).$$ Define $\psi:P'_A \to P'_B$ by %
\begin{align}
\label{E:embedweight} %
\psi(\Lambda_i) = \dfrac{\kappa}{o_i} \,\too_i,
\end{align}
where $\too_i = \sum_{j \in O_i} \omega_j$ and $\tb_i = \sum_{j \in
O_i} \beta_j$ for $i\in I$. We have
\begin{equation} \label{E:embedroot}
\begin{split}
  \psi(\al_i) &= \sum_{i'} a_{i'i} \psi(\La_{i'}) \\
  &= \sum_{i'} a_{i'i} \dfrac{\kappa}{o_{i'}} \sum_{j'\in O_{i'}} \om_{j'} \\
  &= \dfrac{\kappa}{o_i} \sum_{i'} \sum_{j\in O_i} \sum_{j'\in O_{i'}} b_{j'j}
  \om_{j'} \\
  &= \dfrac{\kappa}{o_i} \tb_i
\end{split}
\end{equation}
by \eqref{E:root2weightA}, \eqref{E:embedweight}, and
\eqref{E:asfold}. Define $\varphi:Q^\vee_A\to Q^\vee_B$ by
\begin{align}
\label{E:embedcoroot}%
\varphi(\al_i^\vee) = \tb_i^\vee = \sum_{j\in O_i} \beta_j.
\end{align}

\subsection{Weyl groups}
For $i\in I$ define
\begin{align} \label{E:embedreflect}
  f(s_i^A) = \prod_{j\in O_i} s_j^B \in W_B.
\end{align}
Since $\pi$ is admissible the reflections $\{s_j^B\mid j\in O_i\}$
commute with each other so that the product in
\eqref{E:embedreflect} is independent of the order of its factors.

Steinberg \cite{Stein} showed that there is an embedding of Weyl
groups $W_A\to W_B$. By \cite{Nan} it respects the Bruhat order.

\begin{thm} \label{T:embedWeyl} \
\begin{enumerate}
\item
\cite{Stein} There is an injective group homomorphism $f:W_A\to W_B$
defined by \eqref{E:embedreflect}, whose image is the subgroup
$W_B^\pi$ of $\pi$-fixed elements in $W_B$, where $\pi$ acts on
$W_B$ by $\pi(s_i)=s_{\pi(i)}$.
\item \cite[Prop 3.3 and Thm. 1.2]{Nan} %
$v \le w$ in $W_A$ if and only if $f(v)\le f(w)$ in $W_B$. Moreover,
if $w=s_{i_1}\dotsm s_{i_N}$ is a reduced decomposition in $W_A$
then $f(w)=f(s_{i_1})f(s_{i_2})\dotsm f(s_{i_N})$ is a
length-additive factorization in $W_B$.
\end{enumerate}
\end{thm}

\begin{cor} \label{C:embedstrongcover} \
\begin{enumerate}
\item
$\pi$ acts on $\RealRoot^+(\geh(B))$.
\item
Suppose $v \lessdot w=vs_\al$ in $W_A$ for $\al\in
\RealRoot^+(\geh(A))$. Then there is a unique $\pi$-orbit $O\subset
\RealRoot^+(\geh(B))$ such that the reflections $\{s_\gamma\mid
\gamma\in O\}$ commute, $f(s_\al)=\prod_{\gamma\in O} s_\gamma$, and
there is an isomorphism of the boolean lattice of subsets $O'$ of
$O$ with the interval $[f(v),f(w)]$ in $(W_B,\le)$ given by
$O'\mapsto f(v)\prod_{\gamma\in O'} s_\gamma$. We call $O$ the orbit
associated with the cover $v\lessdot w$.
\item
Let $w\in W_A$ and $\gamma\in\RealRoot^+(\geh(B))$ be such that
$f(w) s_\gamma \lessdot f(w)$ in $W_B$, and let $O\subset
\RealRoot^+(\geh(B))$ be the $\pi$-orbit of $\gamma$. Then there is
a covering relation $v\lessdot w$ of $w$ in $W_A$ of which $O$ is
the associated orbit.
\item
Let $v\in W_A$, $\gamma\in \RealRoot^+(\geh(B))$ be such that $f(v)
\lessdot f(v)s_\gamma$, and $O\subset\RealRoot^+(\geh(B))$ the
$\pi$-orbit of $\gamma$. Then there is a covering relation
$v\lessdot w$ in $W_A$ of which $O$ is the associated cover.
\end{enumerate}
\end{cor}
\begin{proof}
For (1), let $\gamma\in \RealRoot^+(\geh(B))$, with $\gamma=\tu
\beta_j$ for some $\tu\in W_B$ and $j\in J$. Then
$\pi(\gamma)=\pi(\tu) \beta_{\pi(j)}\in \RealRoot(\geh(B))$ and
$\pi$ clearly preserves the set of positive roots, so that $\pi$
acts on $\RealRoot^+(\geh(B))$.

For (2), there is a unique length-additive factorization $w=u_1 s_i
u_2$ in $W_A$ such that $v=u_1u_2$. We have $\al=u\al_i$ and
$s_\al=us_iu^{-1}$ where $u=u_2^{-1}$. Define $O=\{f(u) \beta_j\mid
j\in O_i\}$; since $f(u)$ is $\pi$-invariant, we have
$\pi(f(u)\beta_j)=\pi(f(u))\beta_{\pi(j)}=f(u)\beta_{\pi(j)}\in
O(\alpha)$, so that $O$ is a $\pi$-orbit. For $j\in O_i$ we have
$s_{f(u)\beta_j} = f(u) s_j f(u)^{-1}$, so the reflections
$\{s_\gamma\mid \gamma\in O\}$, commute, being conjugate to
commuting reflections $\{s_j\mid j\in O_i\}$. Since $f$ is a
homomorphism we have $f(s_\al)=f(u) (\prod_{j\in O_i} s_j)
f(u)^{-1}= \prod_{j\in O_i} f(u) s_j f(u)^{-1}=\prod_{\gamma\in
O(\alpha)} s_\gamma$.

By Theorem \ref{T:embedWeyl} $f(w)=f(u_1)(\prod_{j\in O_i}
s_j)f(u_2)$ and $f(v)=f(u_1)f(u_2)$ are length-additive
factorizations. It follows that there is an isomorphism of the
boolean lattice of subsets $S\subset O_i$ with $[f(v),f(w)]$ where
$S\mapsto f(u_1) (\prod_{j\in S} s_j ) f(u_2)$. The desired
isomorphism is given by sending $S\to O'$ where
$O'=\{f(u)\beta_j\mid j\in S\}$.

For (3), let $w=s_{i_1}\dotsm s_{i_k}$ be a reduced decomposition.
Then $f(w)=f(s_{i_1})\dotsm f(s_{i_k})$ is length-additive by
Theorem \ref{T:embedWeyl}. Therefore the cover $f(w)s_\gamma
\lessdot f(w)$ is obtained by removing some unique reflection in
$f(s_{i_r})$ for some unique $r$. Let $u_1=s_{i_1}\dotsm
s_{i_{r-1}}$, $i=i_r$, and $u^{-1}=u_2=s_{i_{r+1}}\dotsm s_{i_k}$.
Letting $\al=u\al_i\in\RealRoot^+(\geh(A))$ we find that the
$\pi$-orbit $O$ of $\gamma$ is the orbit associated with the cover
$v\lessdot w$ where $v=u_1u_2$.

The proof of (4) is similar.
\end{proof}

The following result is proved similarly.

\begin{cor} \label{C:embedweakcover}
Let $v \wless s_i v$ in $W_A$ for some $i\in I$. Then there is a
poset isomorphism from the boolean lattice of subsets $O'$ of $O_i$,
to the interval $[f(v),f(s_i v)]$ of $(W_B,\preceq)$ given by
$O'\mapsto (\prod_{j\in O'} s_j) f(v)$.
\end{cor}

\subsection{Pairings}
The action of $W_B$ on $P_B$ descends to $P'_B$, and similarly for
$W_A$.

\begin{thm} \label{T:embed} \
\begin{enumerate}
\item
For all $\al^\vee\in Q^\vee_A$ and $\la\in P'_A$ we have
\begin{align}
\label{E:dilate}%
\inner{\varphi(\al^\vee)}{\psi(\la)} = \kappa\,
\inner{\al^\vee}{\la}.
\end{align}
\item
For all $\La\in P'_A$ and $w\in W_A$ we have
\begin{equation} \label{E:Weylequiv}
  \psi(w \Lambda) = f(w) \psi(\Lambda).
\end{equation}
\item
Let $w \in W_A$, $\La \in P'_A$ and $i \in I$.  Then
\begin{equation} \label{E:dilateWeyl}
 \inner{\varphi(\al^\vee)}{f(w) \psi(\La)} = \kappa\inner{\al^\vee}{w\La}
\end{equation}
and in particular we have
\begin{equation} \label{E:embedChevalley}
\inner{\tb_i^\vee}{f(w) \psi(\La)} = \kappa
\inner{\al_i^\vee}{w\La}.
\end{equation}
\item The map $\varphi$ sends $Q^\vee_A \cap Z(\geh(A))$ into $Z(\geh(B))$.
\end{enumerate}
\end{thm}
\begin{proof} By linearity it suffices to check \eqref{E:dilate} for
$\al^\vee=\al^\vee_i$ and $\la=\La_k$ for $i,k\in I$. We have
\begin{align*}
  \inner{\varphi(\al^\vee_i)}{\psi(\La_k)} =
  \dfrac{\kappa}{o_k} \inner{\tb_i^\vee}{\too_k} = \dfrac{\kappa}{o_k} o_i \delta_{ik} = \kappa \delta_{ik}
  = \kappa \inner{\al_i^\vee}{\La_k}.
\end{align*}
This implies \eqref{E:dilate}.

It suffices to prove \eqref{E:Weylequiv} for $w=s_{i'}^A$ and
$\La=\La_i$. For $i=i'$ we have
$$
\psi(s_i\Lambda_i) = \psi(\Lambda_i - \alpha_i) =
\dfrac{\kappa}{o_i}(\too_i - \tb_i)
$$
by \eqref{E:embedweight} and \eqref{E:embedroot}. Since $\pi$ is
admissible, $(\prod_{j \in O_i}s_j^B) \too_i = \too_i - \tb_i$,
giving the required result. The case that $i\ne i'$ is even easier.

Equation \eqref{E:dilateWeyl} follows immediately from
\eqref{E:Weylequiv} and \eqref{E:dilate} with $\la=w\La$. Equation
\eqref{E:embedChevalley} is implied by \eqref{E:dilateWeyl} with
$\al=\al_i$ and \eqref{E:embedcoroot}.

Let $K\in Q_A^\vee$. For any $i\in I$ and $j\in O_i$, by
\eqref{E:dilate}, \eqref{E:embedroot}, and the $\pi$-invariance of
$\varphi(K)$, we have
\begin{equation} \label{E:Kcoefs}
  \kappa \inner{\varphi(K)}{\beta_j} =\dfrac{\kappa}{o_i} \inner{\varphi(K)}{\tb_i} = \kappa\inner{K}{\al_i}.
\end{equation}
It follows that $\varphi$ sends $Z(\geh(A))$ into $Z(\geh(B))$.
\end{proof}

\begin{remark} Let $A$ be a GCM of affine type, obtained as in
\cite{L} by folding $(B,\pi)$. Then the canonical central elements
and null roots are related by $\varphi(K^A)=K^B$ and
$\psi(\delta^A)=r^\vee\delta^B$ where $r^\vee$ is the ``twist" of
the dual affine root system $X^{(r^\vee)}_N$ to that of $A$ in the
nomenclature of \cite{Kac}.
\end{remark}

\subsection{Folding and insertion}
\label{SS:foldins}

Let $A=(a_{ij}\mid i,j\in I)$ be a GCM with associated Kac-Moody
algebra $\geh(A)$, $i'\in I$, and $K\in Z^+(\geh(A))$. Suppose $A$
is obtained by folding the GCM $B=(b_{ij}\mid i,j\in J)$ with
admissible automorphism $\pi$. Choose $j'\in O_{i'}$. We shall
construct the dual graded graph $(\Gs^A(\La_{i'}),\Gw^A(K))$ from
the dual graded graph $(\Gs^B(\om_{j'}),\Gw^B(\varphi(K)))$. The
construction only requires the subset $W_B^\pi\subset W_B$ of
$\pi$-invariant vertices, and the edges incident to them, grouped
according to their $\pi$-orbits.

\begin{remark} \label{R:orbitrepchoice}
The choice of $j'\in O_{i'}$ is immaterial; if one chooses another
element of $O_{i'}$ then the resulting type $B$ structures are
transported to each other by a power of the automorphism $\pi$.
\end{remark}

\begin{prop} \label{P:foldstrongedges}
Let $v\lessdot w$ in $W_A$ with $\al\in \RealRoot^+(\geh(A))$ such
that $w=vs_\al$ and let $O\subset\RealRoot^+(\geh(B))$ be the
associated orbit of the cover $v\lessdot w$, defined in Corollary
\ref{C:embedstrongcover}. Then the following sets have the same
cardinality.
\begin{enumerate}
\item
Marked edges $v\overset{m}\to w$ in $\Gs^A(\La_{i'})$.
\item
The disjoint union over $\gamma\in O$ of the sets $\{f(w)s_\gamma
\overset{M}\to f(w)\}$ of marked edges going into $f(w)$ in the
interval $[f(v),f(w)]$ in $\Gs^B(\om_{j'})$.
\item
The disjoint union over $\gamma\in O$ of the sets
$\{f(v)\overset{M}\to f(v) s_\gamma \}$ of marked edges coming out
of $f(v)$ in the interval $[f(v),f(w)]$ in $\Gs^B(\om_{j'})$.
\end{enumerate}
\end{prop}
\begin{proof} By Corollary \ref{C:embedstrongcover} the last two
sets are in bijection. Let $u\in W_A$ and $i\in I$ be such that
$\al=u\al_i$. Then $O=\{f(u)\beta_j\mid j\in O_i\}$. Using Theorem
\ref{T:embed} and the $\pi$-invariance of $\sum_{\gamma\in O}
\gamma^\vee$, we have
\begin{equation} \label{E:embedm}
\begin{split}
\sum_{\gamma\in O} m^B_{\om_{j'}}(f(w)s_\gamma,f(w)) &=
\sum_{\gamma\in O} \inner{\gamma^\vee}{ \om_{j'}} =
\dfrac{1}{o_{i'}} \sum_{\gamma\in O} \inner{\gamma^\vee}{\omt_{i'}} \\
&= \dfrac{1}{\kappa} \sum_{\gamma\in O}
\inner{\gamma^\vee}{\psi(\La_{i'})} \\ &= \dfrac{1}{\kappa}
\sum_{j\in O_i}
\inner{f(u)\beta_j^\vee}{\psi(\La_{i'})} \\
&= \dfrac{1}{\kappa} \inner{\tb_i^\vee}{f(u)^{-1} \psi(\La_{i'})} \\
&= \inner{\al_i^\vee}{u^{-1} \La_{i'}} =
\inner{\al^\vee}{\La_{i'}}=m^A_{\La_{i'}}(ws_\al,w).
\end{split}
\end{equation}
\end{proof}

By definition the graded graph $\Gs^B(\om_{j'})^\pi$ has vertex set
$f(W_A)=W_B^\pi\subset W_B$, grading function $h(f(v))=\ell(v)$ for
$v\in W_A$, and for every cover $v\lessdot w$ in $W_A$, an edge from
$f(v)$ to $f(w)$ whose multiplicity is the common number in
Proposition \ref{P:foldstrongedges}. It is completely specified by
the $\pi$-invariant elements of $\Gs^B(\om_{j'})$ and their incident
edges.

\begin{cor} \label{C:foldstrongiso} The graded graphs
$\Gs^A(\La_{i'})$ and $\Gs^B(\om_{j'})^\pi$ are isomorphic.
\end{cor}
\begin{proof}
By Theorem \ref{T:embedWeyl} the map $f:W_A\to W_B^\pi$ is a grade-
and edge-preserving bijection. The edge multiplicities agree by
Proposition \ref{P:foldstrongedges}.
\end{proof}

We call $\Gs^B(\om_{j'})^\tau$ the \textit{folded strong graph} and
its tableaux \textit{folded strong tableaux}.

\begin{prop} \label{P:foldweakedges}
Let $v\in W_A$ and $i\in I\setminus \Des(v)$ so that $v\prec w=s_i
v$. Fix any $j\in O_i$. Then the following sets have the same
cardinality.
\begin{enumerate}
\item The marked edges $v\overset{m'}\to w$ in $\Gw^A(K)$.
\item The marked edges
$f(v)\overset{M'}\to s_jf(v) $ in $\Gw^B(\varphi(K))$.
\item The marked edges
$s_jf(w) \overset{M'}\to f(w)$ in $\Gw^B(\varphi(K))$.
\end{enumerate}
\end{prop}
\begin{proof}
By a proof similar to the one for \eqref{E:Kcoefs} and recalling the
definition \eqref{E:weaklabel} we have
\begin{equation} \label{E:embedn}
\begin{split}
  n^A_K(v,s_iv) = \inner{K}{\La_i}
  = \inner{\varphi(K)}{\omega_j}
  = n^B_{\varphi(K)}(f(v),s_j f(v)).
\end{split}
\end{equation}
\end{proof}

By definition the graded graph $\Gw^B(\varphi(K))^\pi$ has vertex
set $f(W_A)=W_B^\pi$, grading function $h(f(w))=\ell(w)$, and for
each $v\in W_A$ and $i\in I\setminus \Des(v)$, an edge from $f(v)$
to $f(s_i v)$ whose multiplicity is the common multiplicity in
Proposition \ref{P:foldweakedges}.

\begin{cor} \label{C:foldweakiso}
The graded graphs $\Gw^A(K)$ and $\Gw^B(\varphi(K))^\pi$ are
isomorphic.
\end{cor}

We call $\Gw^B(\varphi(K))^\pi$ the \textit{folded weak graph} and
its tableaux \textit{folded weak tableaux}.

\begin{cor} \label{C:folddual}
$(\Gs^B(\om_{j'})^\pi,\Gw^B(\varphi(K))^\pi)$ and
$(\Gs^A(\La_{i'}),\Gw^A(K))$ are isomorphic dual graded graphs.
\end{cor}
\begin{proof} This follows immediately from Corollaries
\ref{C:foldstrongiso} and \ref{C:foldweakiso}.
\end{proof}

Let $\db^B$ be a differential bijection for
$(\Gs^B(\om_{j'}),\Gw^B(\varphi(K)))$. We shall use $\db^B$ to
construct a differential bijection for the folded dual graded graph
$(\Gs^B(\om_{j'})^\pi,\Gw^B(\varphi(K))^\pi)$, which, by the
identifications given in Propositions \ref{P:foldstrongedges} and
\ref{P:foldweakedges}, yields a differential bijection $\db^A$ for
$(\Gs^A(\La_{i'}),\Gw^A(K))$.

By Remark \ref{R:KMoffdiag} the off-diagonal part $\db_{vw}^A$ for
$v\ne w$ in $W_A$, has already been specified.

For the diagonal terms, since $\Des(f(v))=\bigsqcup_{i\in \Des(v)}
O_i$, by \eqref{E:embedm} in the special case of a cover of the form
$s_i v = vs_\al \gtrdot v$ with $\al=v^{-1}\al_i\in\RealRoot^+$ and
\eqref{E:embedn} we have
\begin{equation} \label{E:DUfold}
\begin{split}
  &\sum_{i\in I\setminus\Des(v)}   m^A_{\La_{i'}}(v,s_iv)\,
    n^A_K(v,s_iv) \\
= &\sum_{j\in J\setminus\Des(f(v))} m^B_{\omega_{j'}}(f(v),s_j
f(v))\, n^B_{\varphi(K)}(f(v),s_j f(v)).
\end{split}
\end{equation}
For $i\in \Des(v)$ and $j\in O_i$, we have $vs_\al =s_i v \prec v$
where $\al=v^{-1}\alpha_i\in -\RealRoot^+$, and an analogous
computation yields
\begin{equation} \label{E:UDfold}
\begin{split}
  &\sum_{i\in \Des(v)}   m^A_{\La_{i'}}(s_iv,v)
    n^A_K(s_i v,v) \\
= &\sum_{j\in \Des(f(v))} m^B_{\omega_{j'}}(s_j f(v),f(v))
n^B_{\varphi(K)}(s_j f(v),f(v)).
\end{split}
\end{equation}
Using the bijections of Proposition \ref{P:foldstrongedges} and
\ref{P:foldweakedges}, we obtain bijections $\DU^A_v \to
\DU^B_{f(v)}$ and $\UD^A_v\to \UD^B_{f(v)}$. Under these
identifications we obtain a differential bijection $\db^A$ for
$(\Gs^A(\La_{i'}),\Gw^A(K))$.

\begin{prop} \label{P:foldbij}
Let the GCM $A=(a_{ij}\mid i,j\in I)$ be obtained by folding from
the GCM $B=(b_{ij}\mid i,j\in J)$ with admissible automorphism
$\pi$. Then for any $i'\in I$ and $j'\in O_{i'}$, a differential
bijection $\db^B$ for $(\Gs^B(\omega_{j'}),\Gw^B(\varphi(K))$
restricts to a differential bijection $\db^A$ for the pair of dual
graded graphs $(\Gs^A(\La_{i'}),\Gw^A(K))$.
\end{prop}

We shall give an extensive example in Section \ref{sec:C}.

\begin{remark}
It is possible to axiomatize conditions for an arbitrary pair of
dual graded graphs $(\Gamma,\Gamma')$ and an automorphism $\pi$ of
$(\Gamma,\Gamma')$ to give rise to a folded insertion in the manner
we have described for Kac-Moody graded graphs. The key properties
needed are abstract graph-theoretic formulations of
Theorem~\ref{T:embedWeyl} and Corollary \ref{C:embedstrongcover}.
Since we have no interesting examples which do not come from
Kac-Moody dual graded graphs, we will not make this precise.
\end{remark}

\section{Affine type $C$ combinatorics}
\label{sec:C}

In this section we consider folded insertion for the affine root
system $C_n^{(1)}$. Folding works for the entire Weyl group. However
here we shall restrict our discussion to the explicit description of
folded insertion for the maximal parabolic quotient of the affine
Weyl group given by the dual graded graphs
$(\Gs(\La_{i'}),\Gw(K))^{i'}$ using $2n$-cores, where $i'\in I$ is
any Dynkin node. In the limit $n\to \infty$ one obtains a new
Schensted bijection for every integer $i$, which for $i=0$ coincides
with ``standard" Sagan-Worley insertion into shifted tableaux.

Let $I=\{0,1,\dotsc,n\}$ be the Dynkin node set and $(a_{ij})$ the
GCM, with $a_{ii}=2$ for $i\in I$, $a_{i,i+1}=a_{i+1,i}=-1$ for
$1\le i\le n-2$, $a_{01}=a_{n,n-1}=-1$, $a_{10}=a_{n-1,n}=-2$, and
other entries zero. Using the recipe in Section \ref{sec:KM}, the
Weyl group $W_n$ has generators $s_i$ for $i\in I$ satisfying
$s_i^2=1$ for $i\in I$ and $(s_i s_j)^{m_{ij}}=1$ for $i,j\in I$
with $i\ne j$, where $m_{01}=m_{n-1,n}=4$, $m_{i,i+1}=3$ for $1\le
i\le n-2$, and $m_{ij}=0$ for $|i-j| \ge 2$.

\subsection{Folding for $C_n^{(1)}$}
Let $A = C_n^{(1)}$ and $B = A_{2n-1}^{(1)}$ denote the two GCMs. We
use the notation of Section~\ref{sec:fold}.

Let $\pi$ be the admissible automorphism of $B$ given by $j \mapsto
2n-j$ where indices are taken modulo $2n$. We index the $\pi$-orbits
by $O_0=\{0\}$, $O_n=\{n\}$, and $O_i=\{i,2n-i\}$ for $i\in
I\setminus\{0,n\}$. It is easy to check that $A$ is obtained from
$(B,\pi)$ by folding.

Let $K$ be the canonical central element for $C_n^{(1)}$. Let $i'\in
I$ and $j'\in O_{i'}$. We define folded insertion for the dual
graded graph $(\Gs^A(\La_{i'}),\Gw^A(K))^{i'}$, realized by LLMS
insertion for $(\Gs^B(\om_{j'}),\Gw^B(\varphi(K)))$.  We call this
induced folded insertion the ``LLMS insertion for $W_n^{i'}$'' (even
though it also depends on $j'$).

\subsection{$2n$-cores}
As before, fix $i'\in I$ and $j'\in O_{i'}$. The elements of the
parabolic quotient $W_n^{i'}$ may be realized by $2n$-cores as
follows.

By Proposition \ref{prop:LLMS} there is a bijection $c:\ts^0\to
\core$. Using a rotational automorphism of the Dynkin diagram of
type $A_{2n-1}^{(1)}$, for any $k\in J$ one may define the
$k$-action of $\ts$ on $\core$, denoted $w\cdot_k \la$, which is the
same as before except that the diagonal of the cell $(i,j)$ is
$j-i+k$. Since the stabilizer of $\vn$ under the $k$-action of $\ts$
on $\core$ is $(\tilde{S}_{2n})_{J\setminus\{k\}}$, there is a
bijection $c_k:\ts^k\to \ts/ (\tilde{S}_{2n})_{J\setminus\{k\}}\to
\core$ defined by $c_k(w)=w\cdot_k \vn$.

Define the map $sc_{i'}: W_n\to \core$ by $w\mapsto f(w)\cdot_{j'}
\vn$, where $f: W_n \to \ts$ is the Weyl group homomorphism of
Section \ref{sec:fold}.  Note that $f(W_n^{i'}) \subset \ts^{J
\setminus O_{i'}}$.

Denote by $\core^{j'}$ the image of $sc_{i'}$.

The following result is the $C_n^{(1)}$-analogue of (part of)
Proposition~\ref{prop:LLMS}.

\begin{prop}\label{P:LLMSC}
The map $sc_{i'}$ restricts to a bijection $W_n^{i'} \to
\core^{j'}$. For $v,w\in W_n^{i'}$ we have $v\le w$ if and only if
$sc_{i'}(v)\subset sc_{i'}(w)$.
\end{prop}
\begin{proof}
The stabilizer of $W=W_n$ acting on $\vn$ is equal to $W_{I -
\{i'\}}$, so the first statement is immediate.  Let $v,w\in
W_n^{i'}$. The following are equivalent: (1) $v\le w$; (2) $f(v) \le
f(w)$; (3) $f(v) (\ts)_{J\setminus\{j\}} \le f(w)
(\ts)_{J\setminus\{j\}}$ for all $j\in O_{i'}$; (4) $f(v)
(\ts)_{J\setminus\{j'\}} \le f(w) (\ts)_{J\setminus\{j'\}}$; (5)
$sc_{i'}(v) \subset sc_{i'}(w)$. (1) and (2) are equivalent by
Theorem \ref{T:embedWeyl}. (2) and (3) are equivalent by Proposition
\ref{P:Bruhatcrit} below applied to the data $f(v)$, $f(w)$, and
$J\setminus O_{i'}$ in $\ts$. Since $f(v)$ and $f(w)$ are
$\pi$-invariant, (3) and (4) are equivalent, because $j'\in O_{i'}$
and the condition for $j$ is invariant as $j$ runs over a
$\pi$-orbit. (4) and (5) are equivalent by Proposition
\ref{prop:LLMS}.
\end{proof}

For a Coxeter group $W$ and a parabolic subgroup $W_J$, the strong
(Bruhat) order denoted ``$\le$''on the quotient $W/W_J$ is the
partial order naturally induced from the strong order on $W^J$. The
following result is due to Deodhar \cite{Deo}.

\begin{prop} \label{P:Bruhatcrit}
Let $W$ be a Coxeter group with simple generators indexed by $P$ and
let $Q \subset P$.  Suppose $x,y \in W^Q$.  Then $x \le y$ if and
only if $xW_{Q'} \le yW_{Q'}$ for every maximal parabolic subgroup
$W_{Q'} \supset W_Q$.
\end{prop}

We describe $\core^{0}$ explicitly in Section \ref{S:SW}, together
with an explicit description of LLMS insertion for $W_n^0$.  It
would be interesting to obtain an explicit description of
$\core^{j'}$ for arbitrary $j'$. The explicit description of the
Chevalley coefficients in a manner similar to
Proposition~\ref{prop:LLMS}, and of LLMS insertion for $W_n^{i'}$
appears to be rather subtle.

\subsection{Large rank limit of folded LLMS insertion}
We now consider the limit of LLMS insertion for $W_n^{i'}$ as $n$
goes to $\infty$, in such a way that the nodes near $0$ in
$A_{2n-1}^{(1)}$ are stable; for this purpose we label these nodes
$\dotsc,-2,-1,0,1,2,\dotsc$.

Let $A_{\pm\infty}$ be the Kac-Moody algebra\footnote{We allow
infinite Dynkin diagrams in a formal manner.} whose Dynkin diagram
has vertex set $J_\infty=\Z$, with Cartan matrix $(b_{ij})$ such
that $b_{ii}=2$ and $b_{i,i+1}=b_{i+1,i}=-1$, and $b_{ij}=0$
otherwise. Let $S_{\pm\infty}$ be its Weyl group: it has generators
$s_j$ for $j\in J_\infty$, with relations $s_j^2=1$,
$(s_js_{j+1})^3=1$, and $(s_is_j)^2=1$ for $|i-j|\ge2$. Then
$S_{\pm\infty}$ acts on partitions: $s_j\cdot \la$ is obtained from
$\la$ by adding the unique $\la$-addable cell in diagonal $j$ if it
exists, and removing the unique $\la$-removable cell in diagonal $j$
if it exists (remembering the shift in diagonal index by $j'$). Then
$S_{\pm\infty}\vn=\Y$ is the set of all partitions and there is a
bijection $c_{j'}:S_{\pm\infty}^{j'}\cong \Y$.

Let $C_\infty$ be the Kac-Moody algebra with Dynkin node set
$I_\infty=\Z_{\ge0}$ and Cartan matrix $a_{ij}$ with $a_{ii}=2$ for
$i\in I_\infty$, $a_{i,i+1}=a_{i+1,i}=-1$ for $i \in
I_\infty\setminus\{0\}$, $a_{01}=-1$ and $a_{10}=-2$. Then its Weyl
group $W_\infty$ has generators $s_i$ for $i\in I_\infty$ with
relations $s_i^2=1$, $(s_0s_1)^4=1$, $(s_is_{i+1})^3=1$ for $i\in
I_\infty\setminus\{0\}$, and $(s_is_j)^2=1$ for $|i-j|\ge2$. As
before, there is an injective homomorphism $f:W_\infty\to
S_{\pm\infty}$ given by $f(s_0)=s_0$ and $f(s_i)=s_is_{-i}$ for
$i>0$.

Then $W_\infty$ acts on partitions via $f$. Define
$\Y^{i'}=W_\infty\cdot\vn$.  The limit of Proposition \ref{P:LLMSC}
gives a bijection $sc_{i'}:W_\infty^{i'}\cong \Y^{i'}$. The strong
and weak orders on $\core\cong \ts^{j'}$ both converge to Young's
lattice $\Y$. The weak order on $\core^{j'} \cong W_n^{i'}$
converges to the weak order on $\Y^{i'}\cong W_\infty^{i'}$, in
which a cover $\la \subset s_i\la$ adds cells in diagonals $i$ and
$-i$ if $i>0$ or just the cell in diagonal $0$ if $i=0$.

\begin{prop} \label{P:components}
Suppose $\la \subset \mu$ is a strong cover in $\Y^{i'} \cong
W_\infty^{i'}$ with $sc_{i'}^{-1}(\la) = w \lessdot w s_\alpha =
sc_{i'}^{-1}(\mu)$. If $s_\alpha$ is conjugate to $s_0$ then the
Chevalley coefficient $\inner{\alpha^\vee}{\La_{i'}}$ is equal to 1
and $\mu/\la$ has a single connected component which is necessarily
a ribbon.  Otherwise suppose $f(s_\alpha) = s_\beta s_{\beta'}$.
Then the Chevalley coefficient $\inner{\alpha^\vee}{\La_{i'}}$ is
equal 1 or 2 depending on whether one or both of $sc(ws_\beta)$ and
$sc(ws_{\beta'})$ strictly contain $\la$.  Furthermore, each strict
containment has a single connected component which is a ribbon.
\end{prop}
\begin{proof}
That the skew partitions in question contain a single connected
component which equals a ribbon follows from the fact that they are
obtained by the action of a reflection on a partition, which always
changes the shape by a ribbon.  This follows from the ``edge
sequence'' discussion in Section \ref{sec:LLMS}.

Suppose $s_\alpha$ is conjugate to $s_0$ and $f(s_\alpha) =
s_\beta$.  Every strong cover in $\Y$ has Chevalley coefficient
equal to 1 or 0 (since a single box is added).  Now write $f(w) =
xy$ so that $f(ws_\alpha) =x s_0 y$ in $S_{\pm \infty}$. The
Chevalley coefficient $\inner{\alpha^\vee}{\La_{i'}}$ is equal to
the Chevalley coefficient $\inner{\beta^\vee}{\om_{j'}}$ which is
equal to the Chevalley coefficient of the cover $y \lessdot s_0y$ in
$S_{\pm \infty}$ (with respect to $\om_{j'}$).  Since $\mu \ne \la$,
we must have $y \cdot \vn \subsetneq (s_0y) \cdot \vn$. Thus the
required Chevalley coefficient must be non-zero, and hence equal to
1.  We have used the calculation $\inner{y^{-1}
\beta_0^\vee}{\om_{j'}} = \inner{\beta_0^\vee}{y \om_{j'}} =
\inner{\beta_0^\vee}{y' \om_{j'}}$, where $y' = c_{j'}^{-1}(y \cdot
\vn)$ is the ``parabolic component'' of $y$.

The proof for $f(s_\alpha) = s_\beta s_{\beta'}$ follows in a
similar manner.  The only delicate issue is to show that if
$\inner{\beta^\vee}{\om_{j'}} = 1$ then $\la \subset sc(ws_\beta)$
is a strict inclusion.  But $\inner{\beta^\vee}{\om_{j'}} = 1$
implies that $(s_\beta \cdot \vn) \ne \vn$ so that $(ws_\beta \cdot
\vn) \neq (w \cdot \vn)$.
\end{proof}

\begin{figure}
\begin{align*}
\xymatrix@=4pt{%
{.} && {.} && {.} && {.} && {.}  \\
&&&&&&&{\otimes} \\
{.} && {.} && {.} && {.} \ar[rr] && {\yng(1)} \\
&&&{\otimes} \\
{.} && {.}\ar[rr] && {\yng(1)} && {\yng(1)} \ar[rr] && {\yng(1,1)} \\
&{\otimes} \\
{.} \ar[rr] && {\yng(1)} \ar[rr] && {\yng(1,1)}  && {\yng(1,1)}\ar[rr] &&{\yng(1,1,1)} \\
&&&&&{\otimes} \\
{.}\ar[rr] && {\yng(1)}\ar[rr] && {\yng(1,1)}\ar[rr] && {\yng(1,2)}
\ar[rr] && {\yng(1,1,1,2)}
}%
\end{align*}
\begin{align*}
P = \tableau{4 \\ 4 \\ 2 \\ 1 & 3} \qquad Q = \tableau{4\\3\\2\\1&4}
\end{align*}
\caption{Folded insertion of $\sigma=4213$ for $i'=j'=1$}
\label{F:ins2}
\end{figure}
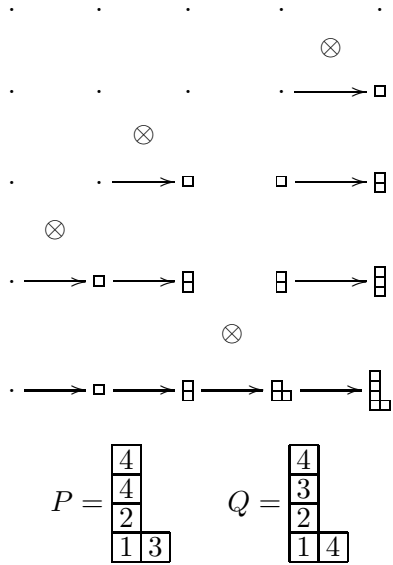

Figure \ref{F:ins2} shows the case of a domino appearing in the
strong tableau $P$, corresponding to the strong cover $s_2s_0s_1
\lessdot s_2s_1s_0s_1$ in $W_\infty^{1}$.

In the case the Chevalley coefficient described in
Proposition~\ref{P:components} is equal to 2, the difference
$\mu/\la$ is a union of two ribbons, since $s_\beta$ and
$s_{\beta'}$ commute.  We have not shown that these ribbons do not
touch, so the difference $\mu/\la$ can potentially be written as the
union of the two ribbons in two ways, corresponding to the left
action of $s_{f(w)\beta} s_{f(w)\beta'}$ and $s_{f(w)\beta'}
s_{f(w)\beta}$. To obtain a strong tableau in $\Y^{i'}$, we must
mark a ribbon for each strong cover which consists of two ribbons.

For the differential bijection, we note that for $\la\in \Y^{i'}$,
once again, $\UD_\la$ is in natural bijection with the set of
$\la$-addable corners and $\DU_\la$ is in natural bijection with the
set of $\la$-removable corners; in this context the corners are
grouped by diagonals of the form $\pm i$ for various $i$. Using the
differential bijection in Example \ref{X:Youngdb}, we obtain a
folded insertion for the limit $\Y^{i'}$ of $\core^{j'}$. It defines
a bijection from permutations $P_n(1)$ to pairs $(P,Q)$ where $P$
and $Q$ are $n$-step strong and weak tableaux with respect to
$\Y^{i'}$.

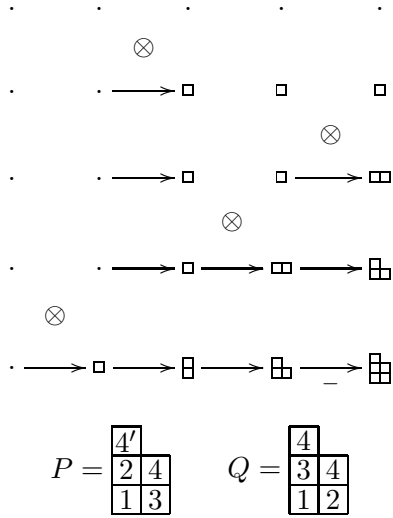
\begin{figure}
\begin{align*}
\xymatrix@=4pt{%
{.} && {.} && {.} && {.} && {.}  \\
&&&{\otimes} \\
{.} && {.} \ar[rr] && {\yng(1)} && {\yng(1)} && {\yng(1)} \\
&&&&&&&{\otimes} \\
{.} && {.}\ar[rr] && {\yng(1)} && {\yng(1)}\ar[rr] && {\yng(2)} \\
&&&&&{\otimes} \\
{.} && {.} \ar[rr] && {\yng(1)} \ar[rr] && {\yng(2)}\ar[rr] &&{\yng(1,2)} \\
&{\otimes} \\
{.}\ar[rr] && {\yng(1)}\ar[rr] && {\yng(1,1)}\ar[rr] && {\yng(1,2)}
\ar[rr]_{-} &&{\yng(1,2,2)}
}%
\end{align*}
\begin{align*}
P = \tableau{4' \\ 2&4\\ 1&3} \qquad Q = \tableau{4\\ 3&4\\ 1&2}
\end{align*}
\caption{Folded insertion of $\sigma=2431$ for $i'=j'=1$}
\label{F:ins1}
\end{figure}

For  example, we let $i'=j'=1$ and compute the folded insertion of
$\sigma=2431$; see the graph $G$ in Figure \ref{F:ins1}. For the
meaning of $\otimes$ see equation \eqref{E:fominsquare} and Example
\ref{X:LLMS}. The arrows represent strong covers, and an arrow is
labeled with $-$ if the strong cover adds two nonadjacent cells and
the marked cell is in a more negative diagonal. The unique arrow
labeled with $-$, corresponds to the entry $4'$ in $P$. Strictly
speaking we should mark one ribbon for each number used, but when
there is no choice we have omitted the marking.
\Yboxdim{4pt}%

Again, with $i'=j'=1$ we compute the folded insertion of $4213$; see
the graph $G$ in Figure \ref{F:ins2}. Note that there is a unique
strong cover that is not a weak cover, corresponding to the domino
in $P$ containing 4s.
\Yboxdim{4pt}%

\subsection{Sagan-Worley insertion}\label{S:SW}
We now consider the important case that $i'=0$. We must have $j'=0$.
The following result is straightforward.

\begin{lem} The set $\core^0$ is the subset of $\core$ of elements
fixed by the transpose $\tr$.
\end{lem}

Using the fact that $f(W_n^0)\subset \ts^0$, the next result is an
easy consequence of Propositions \ref{P:foldbij} and
\ref{prop:LLMS}.  A similar statement holds for $i' = n$.

\begin{prop} \label{P:zero}
Suppose $\la \subset \mu$ is a strong cover in $\core^{0}\cong
W_n^{0}$ with $sc_{0}^{-1}(\mu) = sc_{0}^{-1}(\la)s_\alpha$. Then
the Chevalley coefficient $\inner{\alpha^\vee}{\La_{0}}$ is equal to
the total number of components of $\nu/\la$ for a strong cover $\la
\subset \nu$ in $\core\cong \ts^{0}$ satisfying $\la \subset \nu
\subset \mu$. In particular, if $\la\subset\mu$ is also a left weak
cover in $\core^{0}\cong W_n^{0}$ with $\mu=s_i \la$ for $i\in I$,
then $\inner{\alpha^\vee}{\La_{0}}$ is equal to the number of
$\la$-addable corner cells of residue $i$ or $-i$ modulo $2n$.
\end{prop}

Thus for $i' = 0$, a {\it folded strong tableau} of shape
$\la\in\core^{0}$ is a sequence of strong covers in $\core^{0}\cong
W_n^{0}$ from $\vn$ to $\la$, such that every cover has a marked
component. A {\it folded weak tableau} of shape $\la\in\core^{0}$ is
a sequence of weak covers in $\core^{0}\cong W_n^{0}$ going from
$\vn$ to $\la$; no marking is necessary. With these explicit
descriptions we have

\begin{cor} \label{C:Cfold} LLMS insertion induces a bijection from
the set of permutations $P_n(1)$ to pairs $(P,Q)$ of tableaux of the
same shape $\la\in \core^{0}$, where $\ell(sc_{0}^{-1}(\la))=n$, $P$
is a folded strong tableau and $Q$ a folded weak tableau.
\end{cor}

Again we now consider the $n\to \infty$ limit. In this case the
limit of $\core^0$ is the set $\Y^0$ of partitions fixed under the
transpose. The strong and weak orders both converge to the same
order on $\Y^0$. Since added cells are in transpose-symmetric
positions, when marking a strong cover of the form $\la \lessdot s_i
\la$, one must mark either the added cell in diagonal $i$ or $-i$ if
$i>0$.

Folded weak tableaux $Q$ are in obvious bijection with standard
shifted tableaux $Q^*$, given by taking only the part on one side of
the diagonal. A similarly obvious bijection exists from folded
strong tableaux to standard shifted tableaux in which off-diagonal
entries may or may not have a mark; we choose the bijection so that
a mark in a shifted tableau $P^*$ indicates that the corresponding
cell with negative diagonal index is marked in the folded strong
tableau $P$.

Thus we have the correct kinds of tableaux to compare with
Sagan-Worley insertion.

\begin{ex} Let $n=7$ and $\sigma=2673541$. See Example \ref{X:LLMS}
for the way to interpret $\sigma$ and the symbols $\otimes$. We draw
the graph $G_{ij}$ for the folded insertion of $\sigma$. We draw
arrows to represent strong covers and place an asterisk on an edge
if the marked cell is on a negative diagonal. $P$ and $Q$ are the
folded strong and weak tableaux respectively. $P^*$ and $Q^*$ are
the shifted tableaux corresponding to $P$ and $Q$.

\Yboxdim{4pt}%
\begin{figure}
\begin{align*}
\xymatrix@=4pt{%
{.} && {.} && {.} && {.} && {.} && {.} && {.} && {.} \\
&&&{\otimes} \\
{.} && {.}\ar[rr] && {\yng(1)} && {\yng(1)} && {\yng(1)} &&
{\yng(1)} && {\yng(1)} && {\yng(1)} \\
&&&&&&&&&&&{\otimes} \\
{.} && {.} \ar[rr] && {\yng(1)} && {\yng(1)} &&{\yng(1)} &&{\yng(1)}
\ar[rr] && {\yng(1,2)} && {\yng(1,2)} \\
&&&&&&&&&&&&&{\otimes} \\
{.} && {.}\ar[rr] && {\yng(1)} && {\yng(1)} &&{\yng(1)} &&{\yng(1)}
\ar[rr] && {\yng(1,2)}\ar[rr] && {\yng(1,1,3)} \\
&&&&&{\otimes} \\
{.} && {.}\ar[rr] && {\yng(1)}\ar[rr] && {\yng(1,2)} &&{\yng(1,2)}
&&{\yng(1,2)}
\ar[rr] && {\yng(2,2)}\ar[rr] && {\yng(1,2,3)} \\
&&&&&&&&&{\otimes} \\
{.} && {.}\ar[rr] && {\yng(1)}\ar[rr] && {\yng(1,2)}
&&{\yng(1,2)}\ar[rr] &&{\yng(1,1,3)}
\ar[rr] && {\yng(1,2,3)}\ar[rr] && {\yng(2,3,3)} \\
&&&&&&&{\otimes} \\
{.} && {.}\ar[rr] && {\yng(1)}\ar[rr] && {\yng(1,2)}\ar[rr]
&&{\yng(1,1,3)}\ar[rr] &&{\yng(1,2,3)}
\ar[rr]_{*} && {\yng(2,3,3)}\ar[rr] && {\yng(3,3,3)} \\
&{\otimes} \\
{.}\ar[rr] && {\yng(1)}\ar[rr]_{*} && {\yng(1,2)}\ar[rr] &&
{\yng(2,2)}\ar[rr] && {\yng(1,2,3)}\ar[rr]_{*}
&&{\yng(2,3,3)}\ar[rr]_{*} &&{\yng(1,2,3,4)} \ar[rr] &&
{\yng(1,3,3,4)}
}%
\end{align*}
\newcommand{\prb}{2'}
\newcommand{\pre}{5'}
\newcommand{\prf}{6'}
\Yboxdim{10pt}
\begin{align*}
P &= \tableau{ 6' \\  4 & 5' & 7 \\ 2' & 3 & 5 \\ 1 & 2 & 4' & 6} &
Q &= \tableau{ 7 \\ 3&5&6 \\ 2&4&5 \\ 1&2&3& 7 } \\
P^* &= \young(::7,:3\pre,1\prb4\prf) &
Q^* &= \young(::6,:45,1237)
\end{align*}
\caption{Growth diagram illustrating Theorem~\ref{T:SWinsertion}}
\label{F:growth}
\end{figure}
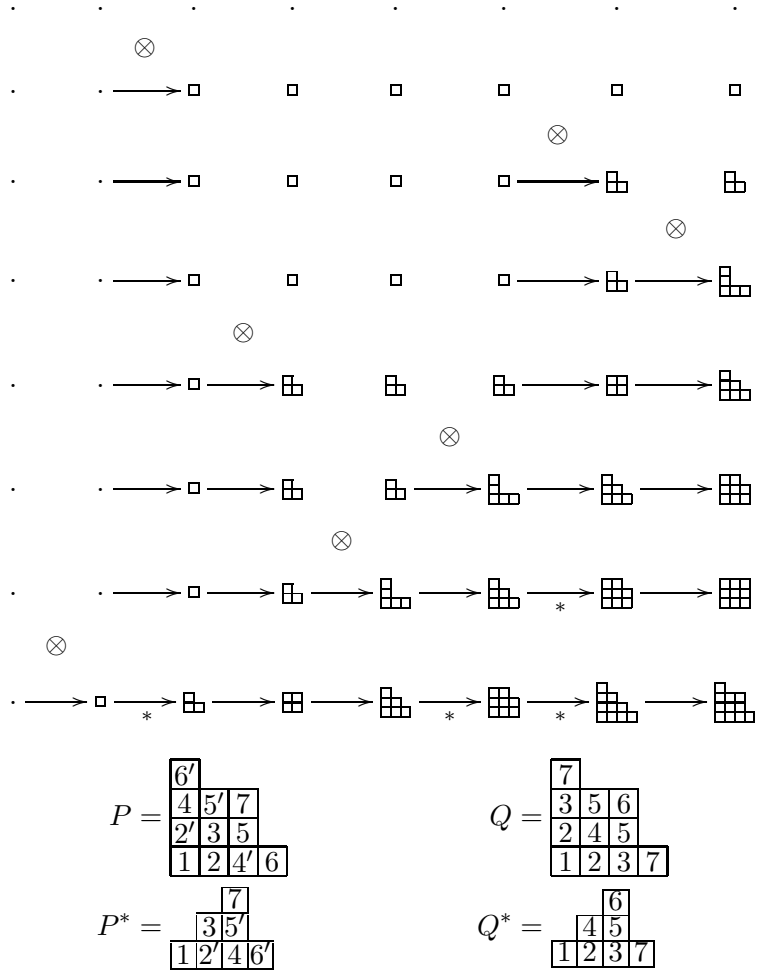
\end{ex}

By reformulating Sagan-Worley insertion using Fomin's setup, we
obtain the following theorem. Note the exchange of $P^*$ and $Q^*$.

\begin{thm} \label{T:SWinsertion}
Let $\sigma$ map to $(P,Q)$ under folded insertion. Then the pair
$(Q^*,P^*)$ of shifted tableaux, is the image of $\sigma^{-1}$ under
Sagan-Worley insertion.
\end{thm}

\begin{remark}
In~\cite[Prop. 6.2]{Ha}, Haiman relates Sagan-Worley shifted
insertion with left-right insertion. It is natural to ask whether
one can connect Sections~\ref{sec:mixed} and~\ref{SS:foldins} in a
similar manner.  Unfortunately, a straightforward generalization of
Haiman's result does not appear to be possible.  For example when
$\pi$ has order 2, one would need to relate the left-right (or
mixed) insertion of a colored permutation on $2r$ letters with
folded insertion of a permutation on $r$ letters.  Length
considerations show that this can be done only if each orbit of
$\pi$ on $J$ has order 2, which nearly never happens in our setup.
\end{remark}

\section{Distributive parabolic quotients}

\subsection{Proctor's classification}
Let $W$ be a finite irreducible Weyl group with simple generators
$\{s_i\mid i\in I\}$ and set of reflections $T$. Recall the
notations $W_J$ and $W^J$ from before Proposition~\ref{P:parabolic}.
We have
$$
W^J = \{w \in W \mid w < ws_i \text{ for any $i \in J$} \}.
$$
In~\cite{Pro}, Proctor classified the cases when $W^J$ is a
distributive lattice under the weak order. In all such cases,
Stembridge \cite{Ste} showed that the weak and strong orders agree
on $W^J$ and that $W_J$ is a maximal parabolic subgroup of $W$, that
is, $J = I \setminus \{i\}$ for some $i\in I$. We call such $W^J :=
W^i$ {\it distributive parabolic quotients}.

\begin{thm}[\cite{Pro}]
\label{thm:Ste} The distributive parabolic quotients are:
\begin{enumerate}
\item
$W \simeq A_n$; $J = I \setminus \{i\}$ for any $i \in I$.
\item
$W \simeq B_n$; $W_J \simeq B_{n-1}$ or $W_J \simeq A_{n-1}$.
\item
$W \simeq D_n$; $W_J \simeq D_{n-1}$ or $W_J \simeq A_{n-1}$.
\item
$W \simeq G_2$; $J = I \setminus \{i\}$ for any $i\in I$.
\end{enumerate}
\end{thm}
In~\cite{Ste}, it is shown that these cases are also exactly the
parabolic quotients $W^J$ of Weyl groups such that every element $w
\in W^J$ is fully commutative, that is, every two reduced
decompositions of $w$ can be obtained from each other using just the
relations of the form $s_is_j=s_js_i$ for $i,j\in I$.

\subsection{Distributive labeled posets}
\label{sec:dist} We need a slightly more precise form of the results
of Proctor and Stembridge. If $Q$ is a finite poset we let $J(Q)$
denote the poset of (lower) order ideals of $Q$.  The poset $J(Q)$
is a distributive lattice and the fundamental theorem of finite
distributive posets \cite{EC2} says that the correspondence $Q
\mapsto J(Q)$ is a bijection between finite posets and finite
distributive lattices. Suppose $P$ is a finite poset and $\omega:
\{x \lessdot y\} \to A$ is a labeling of the edges of the Hasse
diagram of $P$ with elements of some set $A$. We call $(P,\omega)$
an edge-labeled poset. We say that $(P, \omega)$ is a {\it
distributively labeled lattice} if
\begin{enumerate}
\item
$P = J(Q)$ is  a distributive lattice; and
\item
there is a vertex (element) labeling $\pi:Q \to A$ such that
$$
\omega(I\setminus\{q\} \lessdot I) = \pi(q)
$$
for any $I \in J(Q)$ and $q$ maximal in $I$.
\end{enumerate}

If $W$ is a Weyl group, we may label the edges of the Hasse diagram
of the weak order $(W, \prec)$ with simple reflections: the cover $w
\prec s_i w$ is labeled with $s_i$.  We denote the resulting
edge-labeled poset by $W_\weak$.  Similarly define $W_\strong$ to be
the strong order where $w \lessdot wt$ is labeled with $t \in T$.
These labeled posets restrict to give labeled posets $W^i_\weak$ and
$W^i_\strong$.  Note that each cover relation in $W^i$ under either
order is itself a cover relation in $W$. Thus $W^i_\weak$ and
$W^i_\strong$ are induced subgraphs of $W_\weak$ and $W_\strong$.

\begin{thm}
\label{thm:strongSte} Suppose $W^i$ is a distributive parabolic
quotient. Then the strong and weak orders on $W^i$ coincide. In
particular $W^i_\weak$ and $W^i_\strong$ are distributively labeled
lattices.
\end{thm}

Stembridge \cite[Theorem 2.2]{Ste} proved that $W^i_\weak$ is a
distributively labeled lattice. For the sake of completeness we give
a self-contained proof of Theorem \ref{thm:strongSte}.

\subsection{Cominuscule parabolic quotients}
Let $\Phi$ be an irreducible finite root system and $W$ be its Weyl
group.   Let $\Phi = \Phi^+ \sqcup \Phi^-$ denote the decomposition
of the roots into the disjoint subsets of positive and negative
roots.  Let $\theta = \sum_{i \in I} a_i \alpha_i$ denote the
highest root of $\Phi$.  We say that $i \in I$ is {\it cominuscule}
if $a_i = 1$.

It can be checked case-by-case using Theorem~\ref{thm:Ste} that the
distributive parabolic quotients $W^i$ correspond to cominuscule
nodes $i \in I$ except in the cases $W = G_2$, $W^i = B_n/A_{n-1}$
or $W^i = C_n/C_{n-1}$.  In the latter two cases, one may use the
isomorphic quotients given by their duals $C_n/A_{n-1}$ and
$B_n/B_{n-1}$, which are cominuscule.

For now we suppose that a cominuscule node $i \in I$ has been fixed.
If $\alpha$ and $\beta$ are two roots, we say $\alpha \ge \beta$ if
$\alpha-\beta$ is a sum of positive roots. Recall that $\theta$ is
the unique maximal root under this order. Let $\Phi^{(i)}$ denote
the poset of positive roots which lie above $\alpha_i$.  Clearly
$\theta \in \Phi^{(i)}$. The inversion set of $w\in W$ is defined by
$$
\Inv(w) = \{ \alpha \in \Phi^+ \mid w \alpha < 0\}.
$$

\begin{lem}
\label{lem:comin} Suppose $\alpha,\beta\in\Phi^{(i)}$.  Write
$s_\alpha\beta=\beta+k\alpha$ where $k=-\inner{\alpha^\vee}{\beta}$.
Then $k\in \{0,-1,-2\}$ and
\begin{enumerate}
\item If $\alpha$ and $\beta$ are incomparable then
$s_\alpha\beta=\beta$.
\item If $\alpha > \beta$ then $s_\alpha \beta$ is equal to one of the
following: (i) $\beta$; (ii) $-\gamma$ where $\gamma \in
\Phi^+\setminus \Phi^{(i)}$; or (iii) $-\gamma$ where $\gamma
> \alpha$.
\end{enumerate}
\end{lem}
\begin{proof} To obtain the bounds on $k$ we observe that
for all roots $\gamma\in\Phi$, $-\theta\le \gamma\le \theta$, so
that the coefficient of $\alpha_i$ in $\gamma$, lies between the
corresponding coefficients in $-\theta$ and $\theta$, which are $-1$
and $1$ by the assumption that $i$ is cominuscule.

Suppose that $\alpha$ and $\beta$ are incomparable. Then $\beta -
\alpha$ is neither positive nor negative and hence not a root. Since
the roots in $\Phi$ occur in strings, we must have $k=0$.

If $\alpha > \beta$, the three cases correspond to $k = 0$, $k=-1$,
and $k = -2$.
\end{proof}

For our results on distributive parabolic quotients, we require the
following result, which is a slight strengthening of \cite[Prop.
2.1, Lemma 2.2]{TY}. We include a self-contained proof, part of
which is the same as the proof of \cite[Prop. 2.1]{TY}. In
particular we prove directly that the edge labeled poset
$W^i_\strong$ defined in section \ref{sec:dist} is a distributively
labeled lattice.

\begin{prop}
\label{prop:bij} The map $w \longmapsto \Inv(w)$ defines an
isomorphism of posets $\Invj: (W^i,\le) \to J(\Phi^{(i)})$.
Moreover, if $u\lessdot w$ for $u,w\in W^i$, then writing
$w=us_\alpha$ for $\alpha\in\Phi^+$, we have $\alpha\in\Phi^{(i)}$
and $\Inv(w)=\Inv(u)\sqcup \{\alpha\}$.
\end{prop}
\begin{proof}
Let $w \in W^i$.  First we show that $\Inv(w) \subset \Phi^{(i)}$.
Suppose that $\gamma \in \Inv(w)\setminus\Phi^{(i)}$. If $\gamma =
\alpha_k$ where $k \neq i$ this means $w s_k < w$ which contradicts
the assumption that $w\in W^i$. Otherwise $\gamma = \delta + \rho$
where $\delta,\rho \in \Phi^+\setminus\Phi^{(i)}$. Since $w \gamma <
0$ we have $w \delta < 0$ or $w \rho < 0$ so the same argument
applies. Repeating we obtain a contradiction.

Now we show that $\Inv(w)\in J(\Phi^{(i)})$.  Suppose $\alpha \in
\Inv(w)$ and $\beta < \alpha$.  Then $\gamma = \alpha - \beta\in
\Phi^+\backslash \Phi^{(i)}$ since the coefficient of $\alpha_i$ in
$\gamma$ is zero. Since $\Inv(w) \subset \Phi^{(i)}$, we have
$\gamma\not\in \Inv(w)$, that is, $w\alpha - w\beta = w \gamma
> 0$. Since $w \alpha < 0$ this shows that $w \beta < 0$ as desired.
Thus $\Invj$ is well-defined.

Next we show that $\Invj$ sends covers to covers. Let $u\lessdot w$
with $u,w\in W^i$ and $\alpha\in\Phi^+$ such that $w=us_\alpha$.
Then $0>w\alpha=-u \alpha$ so $\alpha\in \Inv(w)\setminus\Inv(u)$.
For all $\beta\in\Inv(u)$, since $\Inv(u)\in J(\Phi^{(i)})$,
$\alpha>\beta$ or $\alpha>\beta$. Either way we have $w \beta = u
s_\alpha \beta < 0$, since by Lemma~\ref{lem:comin}, $s_\alpha
\beta$ is either equal to $\beta$ or $-\gamma$ for $\gamma \in
\Phi^+\setminus\Inv(u)$. That is, $\Inv(u)\subset\Inv(w)$. Since
$|\Inv(w)|=|\Inv(u)|+1$ it follows that
$\Inv(w)=\Inv(u)\sqcup\{\alpha\}$, so that $\Inv(u)\subset \Inv(w)$
is a covering relation in $J(\Phi^{(i)})$.

Next we show that every covering relation in $J(\Phi^{(i)})$ is the
image of a covering relation in $W^i$, and in particular, that
$\Invj$ is onto. An arbitrary covering relation in $J(\Phi^{(i)})$
is given by $S\setminus\{\alpha\} \subset S$ where $S\in
J(\Phi^{(i)})$ and $\alpha$ is maximal in $S$.

By induction there is a $u\in W^i$ such that
$\Inv(u)=S\setminus\{\alpha\}$. Let $w=us_\alpha$. It suffices to
show that
$$\Inv(w) = S \ \ {\rm and} \ \  w \in W^i.$$
The second claim follows from the first since none of the $\alpha_k$
for $k \neq i$ lie in $\Inv(w)$. For the first claim, since
$\alpha\in\Phi^{(i)}\setminus\Inv(u)$, we may argue as before to
show that $S=\Inv(u)\sqcup \{\alpha\} \subset \Inv(w)$.

For the opposite inclusion, suppose $\beta\in \Phi^+\setminus S$. We
must show that $w\beta > 0$.  Write $s_\alpha\beta=\beta+k\alpha$
for $k\in\Z$. If $k=0$ then we are done as before. If $k>0$ then
$s_\alpha\beta > \alpha$, so that $s_\alpha\beta\in\Phi^+\setminus
S$ since $S$ is an order ideal. But then $s_\alpha\beta\notin
\Inv(u)$ so $w \beta > 0$. So we may assume that $k<0$.

Suppose first that $\beta\in\Phi^{(i)}$. We may assume that $\alpha$
and $\beta$ are comparable by Lemma \ref{lem:comin}. Since $S$ is an
order ideal we have $\beta>\alpha$. If $k=-1$ then
$s_\alpha\beta=\beta-\alpha\in\Phi^+\setminus \Phi^{(i)}$ since the
coefficient of $\alpha_i$ is $1$ in both $\alpha$ and $\beta$. In
particular $s_\alpha\beta\not\in \Inv(u)$ so $w \beta >  0$.
Otherwise $k=-2$. Then $s_\alpha\beta=\beta-2\alpha<0$. We have $0 <
\beta-\alpha < \alpha$ and
$-s_\alpha\beta=2\alpha-\beta=\alpha-(\beta-\alpha)<\alpha$. Since
$S$ is an order ideal it follows that $-s_\alpha\beta\in \Inv(u)$
and $w\beta=us_\alpha\beta>0$ as desired.

Otherwise $\beta\in\Phi^+\setminus\Phi^{(i)}$. Since $i$ is
cominuscule we have $k\in\{-1,0,1\}$. We assume $k=-1$ as the other
cases were already done. Then $s_\alpha\beta=\beta-\alpha<0$ since
its coefficient of $\alpha_i$ is $-1$. Moreover
$\alpha-\beta\in\Phi^{(i)}$. Since $\alpha>\alpha-\beta$ and $S$ is
an order ideal, it follows that $\alpha-\beta\in \Inv(u)$. Therefore
$w\beta=us_\alpha\beta>0$ as desired.

We have shown that every cover in $J(\Phi^{(i)})$ is the image under
$\Invj$ of a cover in $(W^i,\le)$.

The bijectivity of $\Invj$ follows by induction and the explicit
description of the image of a cover under $\Invj$.
\end{proof}

\begin{proof}[Proof of Theorem~\ref{thm:strongSte}]
For the case $W = G_2$, both labeled posets $W^i_\weak$ and
$W^i_\strong$ are chains, so the result follows immediately.  Thus
we may assume that $W^i$ is a cominuscule parabolic quotient.

For $W^i_\strong$ the result follows from
Proposition~\ref{prop:bij}. We label the vertices of $\Phi^{(i)}$ by
reflections, defining $\pi:\Phi^{(i)}\to T$ by $\pi(\alpha) =
s_\alpha$. Each cover $w \lessdot w s_\alpha$ in $W^i_\strong$
corresponds to adding $\alpha \in \Phi^{(i)}$ to $\Inv(w)$.  Thus
the edge label of $w \lessdot w s_\alpha$ agrees with the vertex
label $\pi(\alpha) = s_\alpha$.

For the weak order $W^i_\weak$ let us consider two covers $w
\lessdot w s_\alpha = s_\beta w$ and $v \lessdot v s_\alpha =
s_{\beta'} v$ which have the same label $s_\alpha$ in $W^i_\strong$.
We claim that $s_\beta = s_{\beta'} = s_k$ for some $k \in I$.  The
elements $w$ and $v$ differ by right multiplication by some
$s_\gamma$'s where $\gamma \in \Phi^{(i)}$ is incomparable with
$\alpha$; this is accomplished by passing between $w$ or $v$ to the
element $u\in W^i$ such that $\Inv(u)=\Inv(w)\cap\Inv(v)$. By
Lemma~\ref{lem:comin} these $s_\gamma$'s commute with $s_\alpha$,
and so $w \alpha = v \alpha$. This gives us a map $f: \Phi^{(i)} \to
\Phi^+$ defined by $f(\alpha) = \beta = w\alpha$, which does not
depend on $w \in W^i$ as long as $w \lessdot w s_\alpha$.

To show that $f(\alpha)$ is simple for each $\alpha \in \Phi^{(i)}$,
consider a reduced word $w s_\alpha = s_{k_1} s_{k_2} \cdots
s_{k_l}$. We know that $w^{(r)} = s_{k_r} \cdots s_{k_l} \in W^i$
and that $\Inv(w^{(r)})$ differs from $\Inv(w^{(r+1)})$ by some root
in $\Phi^{(i)}$ since $w^{(r+1)}\lessdot w^{(r)}$. For some value $r
= r^*$, this root is $\alpha$ and by the well-definedness just
proved $f(\alpha) = \alpha_{k_{r^*}}$, since $w^{(r^*)} =
w^{(r^*+1)}s_\alpha$.  This shows that the strong order and weak
order on $W^i$ coincide, and that $W^i_\weak$ is isomorphic to the
poset of order ideals of $\Phi^{(i)}$ where $\Phi^{(i)}$ is labeled
with $\pi(\alpha) = f(\alpha)$.
\end{proof}

\section{Distributive subgraphs of Kac-Moody graded graphs}
In this section we apply Theorem~\ref{thm:strongSte} to the dual
graded graphs constructed in Section~\ref{sec:DGG}.

Let $\geh=\geh(A)$ be the Kac-Moody algebra associated to the
generalized Cartan matrix $A$ and let $W$ be its Weyl group.  Let
$W_\fin \subset W$ be a finite parabolic subgroup corresponding to
some index set $I' \subset I$.  Now suppose that $W_\fin$ has a
distributive parabolic quotient as in Theorem~\ref{thm:Ste}
corresponding to $J = I' \setminus \{i\} \subset I'$.  We let $W^J
\subset W_\fin$ denote the distributive parabolic quotient (we use
$W^J$ instead of $W^i$ in this section since $W^J$ is not a maximal
parabolic quotient of $W$, but of $W_\fin$).

Now let $(\La,K) \in P^+ \times Z^+$ and $(\Gs(\La),\Gw(K))$ be the
pair of dual graded graphs constructed in Section~\ref{sec:DGG}.  By
restricting to the set of vertices $W^J \subset W_\fin \subset W$ we
obtain the induced pair of graded graphs $(\Gs(\La),\Gw(K))^J$.
These graded graphs are not dual (see Remark \ref{R:dualinfinite})
but they still have rich combinatorics.

The distributive lattice $(W^J,\le)$ has two edge labelings. Recall
that in $W^J_\strong$, the edge $v\lessdot w=vs_\al$ is labeled
either by the reflection $s_\al$, while in the strong Kac-Moody
subgraph $\Gs^J(\La)$, the edge $v\lessdot w=vs_\al$ is labeled by
the integer $\inner{\alpha^\vee}{\La}$. Similarly the distributive
lattice $(W^J,\preceq)$ has two edge labelings; in $W^J_\weak$, the
edge $v \prec s_j v$ is labeled by the simple reflection $s_j$,
while in the weak Kac-Moody subgraph $\Gw^J(K)$, the edge $v\prec
s_j v$ is labeled by the integer $\inner{K}{\La_j}$. The following
result is an immediate consequence of Theorem \ref{thm:strongSte}.

\begin{thm}
The induced graded subgraphs $\Gs^J(\La)$ and $\Gw^J(K)$ are
distributively labeled lattices.
\end{thm}

Thus $\Gs^J(\La)$ (resp. $\Gw^J(K)$) can be thought of as the poset
of order ideals in some integer labeled poset $P^J$ (resp. $Q^J$).
The $\La$-strong and $K$-weak tableaux can be thought of as linear
extensions of $P^J$ and $Q^J$ with additional markings.

In the rest of the paper, we give examples of the posets $P^J$ and
$Q^J$ and relate them to classically understood tableaux. In each
case we let $\geh$ be of untwisted affine type,
$I_\fin=I\setminus\{0\}$ and $J=I\setminus\{i\}$ for a fixed node
$i\in I_\fin$ to be specified. We use the canonical central element
$K_\can=\sum_{i\in I} a_i^\vee \alpha_i^\vee$ for $K$ and $\La_i$
for the dominant weight. In this case $P^J$ and $Q^J$ are both
labelings of the poset $\Phi^{(i)}\subset\Phi^+$ for the simple Lie
algebra $\geh_\fin$ whose Dynkin diagram is the subdiagram of that
of $\geh$ given by removing the $0$ node. These examples, with the
exception of $G_2$, can be viewed as providing some additional data
for the posets $\Phi^{(i)}$, whose unlabeled versions were given
explicitly in \cite{TY}. As in \cite{TY} we shall rotate the labeled
Hasse diagrams clockwise by 45 degrees so that the minimal element
is in the southwest corner.  In the following, we let $V^J_\weak,
V^J_\strong$ denote the vertex-labeled posets such that $W^J_\weak =
J(V^J_\weak)$ and $W^J_\strong = J(V^J_{\strong})$.

\begin{center}
\begin{figure}
\begin{tabular}{||c|c||}\hline
\text{Root system} & \text{Dynkin Diagram} \\\hline \hline $A_n$ &
\setlength{\unitlength}{3mm}
\begin{picture}(11,3)
\multiput(0,1.5)(2,0){6}{$\circ$}
\multiput(0.55,1.85)(2,0){5}{\line(1,0){1.55}}
\put(6,1.5){$\bullet$} \put(0,0){$1$} \put(2,0){$2$}
\put(3.5,0){$\cdots$} \put(6,0){$i$} \put(7.5,0){$\cdots$}
\put(10,0){$n$}
\end{picture}
\\ %
\hline $C_n, n\geq 3$ & \setlength{\unitlength}{3mm}
\begin{picture}(11,3) \multiput(0,1.5)(2,0){6}{$\circ$}
\multiput(0.55,1.85)(2,0){4}{\line(1,0){1.55}}
\multiput(8.55,1.75)(0,.2){2}{\line(1,0){1.55}} \put(8.85,1.53){$<$}
\put(10,1.5){$\bullet$} \put(0,0){$1$} \put(2,0){$2$}
\put(4,0){$\cdots$} \put(7,0){$\cdots$} \put(10,0){$n$}
\end{picture}
\\ \hline %
$D_n, n\geq 4$ & $\begin{array}{c} \setlength{\unitlength}{2.9mm}
\setlength{\unitlength}{2.9mm} \begin{picture}(11,3.5)
\multiput(0,1.6)(2,0){5}{$\circ$}
\multiput(0.55,2)(2,0){4}{\line(1,0){1.55}}
\put(8.5,1.95){\line(2,-1){1.55}} \put(8.5,1.95){\line(2,1){1.55}}
\put(10,2.5){$\circ$} \put(10,0.7){$\circ$} \put(0,0){$1$}
\put(2,0){$2$} \put(4,0){$\cdots$} \put(7,0){$\cdots$}
\put(9.1,0){$n\!-\!1$} \put(11, 2.3){$n$} \put(10,2.45){$\circ$}
\put(10,0.75){$\bullet$}
\end{picture}
\end{array}$
\\ \hline
$E_6$ & \setlength{\unitlength}{3mm}
\begin{picture}(9,3.6)
\multiput(0,0.5)(2,0){5}{$\circ$}
\multiput(0.55,0.95)(2,0){4}{\line(1,0){1.6}} \put(0,0.5){$\bullet$}
\put(8,0.5){$\circ$} \put(4,2.6){$\circ$}
\put(4.35,1.2){\line(0,1){1.5}} \put(0,-.6){$1$} \put(2,-0.6){$3$}
\put(4,-.6){$4$} \put(6,-.6){$5$} \put(5,2.5){$2$} \put(8,-.6){$6$}
\end{picture}
\\[1mm] \hline %
$E_7$ & \setlength{\unitlength}{3mm}
\begin{picture}(11,4)
\put(0,0.9){$\circ$} \multiput(2,0.9)(2,0){4}{$\circ$}
\put(10,0.9){$\bullet$}
\multiput(0.55,1.35)(2,0){5}{\line(1,0){1.6}} \put(10,0.9){$\circ$}
\put(4,3){$\circ$} \put(4.35,1.6){\line(0,1){1.5}} \put(0,-.2){$1$}
\put(2,-0.2){$3$} \put(4,-.2){$4$} \put(6,-.2){$5$} \put(5,2.9){$2$}
\put(8,-.2){$6$} \put(10,-.2){$7$}
\end{picture}
\\[1mm] \hline
\end{tabular}
\caption{\label{F:comins} Some cominuscule parabolic quotients}
\end{figure}
\end{center}

\subsection{Type $A_n^{(1)}$} Let $i\in I_\fin$ be
arbitrary. The poset $\Phi^{(i)}$ consists of elements
$\alpha_{p,q}=\alpha_p+\dotsm+\alpha_q$ for $1\le p \le i \le q \le
n$. The weak labeling of $\Phi^{(i)}$ is given by
$\alpha_{p,q}\mapsto s_{p+q-i}$. For example, for $n=7$ and $i=3$
and abbreviating $\alpha_{p,q}$ by $pq$ and $s_j$ by $j$, the
labelings of $\Phi^{(i)}$ by positive roots and simple reflections
are given by:
\begin{align*}
V^J_\strong = \tableau{13&14&15&16&17
\\
23&24&25&26&27
\\
33&34&35&36&37
\\
}\qquad V^J_\weak =
\tableau{1&2&3&4&5 \\
2&3&4&5&6 \\
3&4&5&6&7}
\end{align*}
All labelings in $P^J$ and $Q^J$ are given by the constant $1$. The
resulting strong and weak tableaux are usual standard tableaux.

\subsection{Type $C_n^{(1)}$}
Let $i=n$. Let $\alpha_i=e_i-e_{i+1}$ for $1\le i\le n-1$ and
$\alpha_n=2e_n$ where $e_i$ is the $i$-th standard basis element of
the weight lattice $\Z^n$. Then $\Phi^{(n)}$ consists of the roots
$\alpha_{i,j}=e_i+e_j$ for $1\le i\le j\le n$. We have $a_i^\vee=1$
for all $i$. For $n=4$ we have
\begin{align*}
V^J_\strong = \tableau{14&13&12&11\\ 24&23&22\\ 34&33 \\ 44} %
\quad %
V^J_\weak = \tableau{1&2&3&4\\2&3&4\\3&4\\4} \quad %
P^J = \tableau{2&2&2&1\\2&2&1\\2&1\\1} \quad%
Q^J = \tableau{1&1&1&1\\1&1&1\\1&1\\1}
\end{align*}
The strong tableaux are shifted standard tableaux with two kinds of
markings on offdiagonal entries; these are the standard recording
tableaux for shifted insertion \cite{Sagan}. The weak tableaux are
standard shifted tableaux.

\subsection{Type $D_n^{(1)}$} Let $i=n$. Letting
$\alpha_i=e_i-e_{i+1}$ for $1\le i\le n-1$ and
$\alpha_n=e_{n-1}+e_n$, the roots of $\Phi^{(n)}$ are given by
$\alpha_{p,q}=e_p+e_q$ for $1\le p<q\le n$. We have $a_j^\vee=1$ for
$j\in\{0,1,n-1,n\}$ and $a_j^\vee=2$ otherwise. For $n=5$ we give
the labelings of $\Phi^{(n)}$ below. Note the $1$ in the upper left
corner of $Q^J$.
\begin{align*}
V^J_\strong = \tableau{15&14&13&12\\ 25&24&23\\ 35&34 \\ 45} %
\quad %
V^J_\weak = \tableau{1&2&3&4\\2&3&5\\3&4\\5} \quad %
P^J = \tableau{1&1&1&1\\1&1&1\\1&1\\1} \quad%
Q^J = \tableau{1&2&2&1\\2&2&1\\2&1\\1}
\end{align*}

\subsection{Type $E$}
The computations in this section were made using Stembridge's
Coxeter/Weyl package \cite{Stem:CoxWeyl}. In both of the following
cases, $P^J$ has all labels $1$.

For $E_6^{(1)}$ and $i=1$ with the Dynkin labeling in Figure
\ref{F:comins},
\begin{align*}
V^J_\weak = \tableau{ &&&1&3&4&5&6\\
&&&3&4&2\\
&&2&4&5\\
1&3&4&5&6 }\quad %
Q^J = \tableau{&&&1&2&3&2&1\\ &&&2&3&2\\&&2&3&2\\1&2&3&2&1}.
\end{align*}
For $E_7^{(1)}$ and $i=7$ with the Dynkin labeling in Figure
\ref{F:comins},
\begin{align*}
V^J_\weak =
\tableau{&&&&&&&&7\\&&&&&&&&6\\&&&&&&&&5\\&&&&&&&2&4\\&&&&7&6&5&4&3\\&&&&6&5&4&3&1\\&&&&5&4&2&\\&&&2&4&3\\7&6&5&4&3&1}
\quad
Q^J = %
\tableau{&&&&&&&&1\\&&&&&&&&2\\&&&&&&&&3\\&&&&&&&2&4\\&&&&1&2&3&4&3\\&&&&2&3&4&3&2\\&&&&3&4&2\\&&&2&4&3\\1&2&3&4&3&2}
\end{align*}

\subsection{Type $G_2^{(1)}$}
This case does {\it not} correspond to a cominiscule root.  Pick $i
= 1$ and let $\alpha_1, \alpha_2$ be the two simple roots, so that
the highest root is $3 \alpha_1 + 2\alpha_2$.  Then $a_1^\vee = 1$
and $a_2^\vee = 2$. Abbreviating the reflection $s_{p \alpha_1 + q
\alpha_2}$ by $pq$, we have
\begin{align*}
 V^J_\strong &= \tableau{1& 31 & 21 & 32
& 11} &\quad P^J &= \tableau{1&3&2&3&1} \\
V^J_\weak &= \tableau{ 1&2&1&2&1  } &\quad %
Q^J &= \tableau{ 1&2&1&2&1}
\end{align*}

\end{document}